\documentclass{amsart}
\usepackage{amssymb,latexsym}
\usepackage{amsmath}
\usepackage{graphicx}

\newtheorem{thm}{Theorem}[section]
\newtheorem{cor}[thm]{Corollary}

\newtheorem{defin}[thm]{Definition}
\newtheorem{lemma}[thm]{Lemma}

\newtheorem{prop}[thm]{Proposition}

\newcommand{\aaa}{\mbox{$\alpha$}}

\newcommand{\mcCP}{\mbox{$\mathcal{C}_P$}}
\newcommand{\mcCQ}{\mbox{$\mathcal{C}_Q$}}
\newcommand{\mcB}{\mbox{$\mathcal{B}$}}
\newcommand{\mcC}{\mbox{$\mathcal{C}$}}
\newcommand{\bbb}{\mbox{$\beta$}}

\newcommand{\ddd}{\mbox{$\delta$}}

\newcommand{\gggg}{\mbox{$\gamma$}}
\newcommand{\Ggg}{\mbox{$\Gamma$}}
\newcommand{\kkk}{\mbox{$\kappa$}}

\newcommand{\bdd}{\mbox{$\partial$}}

\newcommand{\Ddd}{\mbox{$\Delta$}}

\newcommand{\xX}{\mbox {\sc x}}
\newcommand{\yY}{\mbox {\sc y}}
\newcommand{\aA}{\mbox {\sc a}}
\newcommand{\bB}{\mbox {\sc b}}

\begin{document}

\subjclass{57M25, 57M27, 57M50}

\keywords {bridge position, Heegaard splitting,
strongly irreducible, weakly incompressible, Conway spheres}

\title{Conway products and links with multiple bridge surfaces}

\author{Martin Scharlemann}
\address{\hskip-\parindent
 Martin Scharlemann\\
  Mathematics Department \\
University of California, Santa Barbara \\
Satna Barbara, CA 93117, USA}
\email{mgscharl@math.ucsb.edu}

\author{Maggy Tomova}
\address{\hskip-\parindent
 Maggy Tomova\\
  Mathematics Department \\
University of Iowa \\
 Iowa city, IA 52240, USA}
\email{mtomova@math.uiowa.edu}

\thanks{Research of the first author is partially supported by an NSF grant.} 

\date{\today}

\begin{abstract}
Suppose a link $K$ in a $3$-manifold $M$ is in bridge position with respect to two different bridge surfaces $P$ and $Q$, both of which are c-weakly incompressible in the complement of $K$.  Then either 
\begin{itemize}
\item $P$ and $Q$ can be properly isotoped to intersect in a nonempty collection of curves
that are essential on both surfaces, or 
\item $K$ is a Conway product with respect to an incompressible Conway sphere that naturally decomposes both $P$ and $Q$ into bridge surfaces for the respective factor link(s). 
\end{itemize}
\end{abstract}

\maketitle

\section{Introduction}

A link $K$ in a 3-manifold $M$ is said to be in {\em bridge position}
with respect to a Heegaard surface $P$ for $M$ if each arc of $K-P$ is
parallel to $P$. $P$ is then called a {\em bridge surface} for $K$ in $M$.  Given a link in bridge position with respect to $P$, it is easy to construct more complex bridge surfaces for $K$ from $P$; for example
by stabilizing the Heegaard surface $P$ or by perturbing $K$ to  introduce a minimum and an adjacent
maximum.  As with Heegaard splitting surfaces for a manifold, it is likely that most links have multiple bridge surfaces even apart from these simple operations.  In an effort to understand how two bridge surfaces for the same link might compare, it seems reasonable to follow the program used in \cite{RS} to compare distinct Heegaard splittings of the same non-Haken $3$-manifold.  The restriction to non-Haken manifolds ensured that the relevant Heegaard splittings were strongly irreducible.   In our context the analogous condition is that the bridge surfaces are c-weakly incompressible (defined below).  The natural analogy to the first step in \cite{RS} would be to demonstrate that any two distinct c-weakly incompressible bridge surfaces for a link $K$ in a closed orientable $3$-manifold $M$ can be isotoped so that their intersection consists of a non-empty collection of curves, each of which is essential (including non-meridinal) on both surfaces. In some sense the similar result in \cite{RS} could then be thought of as the special case in which $K=\emptyset$.  

Here we demonstrate that this is true, so long as there are no incompressible Conway spheres  for the knot $K$ in $M$ (cf Section \ref{sec:Conway} and \cite{GL}). In the presence of Conway spheres a slightly different outcome cannot be ruled out:  the bridge surfaces each intersect a collar of a Conway sphere in a precise way; outside the collar the bridge surfaces intersect only in curves that are essential on both surfaces; and inside the collar there is inevitably a single circle intersection that is essential in one surface and meridional, hence inessential,  in the other.  

\section{Definitions and notation}

Suppose that $K$ is a properly embedded collection of $1$-manifolds in a compact manifold $M$.  For $X$ any subset of $M$, let $X_K$ denote $X-K$. A disk that meets $K$ transversally and at most once is called a {\em (punctured) disk}. Thus the parenthetical (punctured) means {\em either} unpunctured or with a single puncture.   

Suppose $F$ is a compact surface in $M$ transverse to $K$.  An isotopy of $F_K$ will mean an isotopy of $F$ in $M$, fixing $K$ set-wise, so that $F$ is always transverse to $K$.  In particular, it is a proper isotopy of $F_K$ in $M_K$.  A simple closed curve on $F_K$ is {\em essential} if it doesn't bound a (punctured) disk on $F_K$.  (In particular, a closed curve on $F_K$ that bounds a once-punctured disk is considered to be inessential.)  An embedded disk $D \subset M_K$ is a {\em compressing disk} for $F_K$ if $D \cap F_K =\bdd D$ and $\bdd D$ is an essential curve in $F_K$. A {\em cut-disk} for $F_K$ (case (3) of \cite[Definition 2.1]{BS}) is an embedded once-punctured disk $D^c$ in $M_K$ such that $D^c \cap F_K =\bdd D^c$ and $\bdd D^c$ is an essential curve in $F_K$. A (punctured) disk that is either a cut disk or a compressing disk will be called a {\em c-disk} for $F_K$.

Any term describing the compressibility of a surface can be extended to account not only for compressing disks but all c-disks. For example we will call a surface {\em c-incompressible} if it has no c-disks. If $F$ is a splitting surface for $M$ (that is, $M$ is the union of two $3$-manifolds along $F$) we will call $F_K$ {\em c-weakly incompressible} if any pair of c-disks for $F_K$ on opposite sides of the surface intersect. If $F_K$ is not c-weakly incompressible, it is {\em c-strongly compressible}.  (Note the perhaps confusing comparison with standard Heegaard surface terminology:  a Heegaard splitting is weakly reducible (cf \cite{CG})  if and only if the Heegaard surface is strongly compressible.)

A properly embedded arc $(\ddd, \bdd \ddd) \subset (F_K, \bdd F_K)$ is inessential if there is a disk in $F_K$ whose boundary is the end-point union of $\ddd$ and a subarc of $\bdd F$.  Otherwise $\ddd$ is {\em essential}.  The surface $F_K$ is $\bdd$-compressible if there is a disk $D \subset M$ so that the boundary of $D$ is the end-point union of two arcs,  $\ddd = D \cap F_K$ an essential arc in $F_K$ and $\bbb = D \cap \bdd M$. Note that arcs in $F$ with one or more end-points at $K \cap F$ do not arise in either context, though such arcs will play a role in our argument.  If $F$ is a properly embedded twice-punctured disk in $M$ then a $\bdd$-compression of $F_K$ may create two $c$-compressing disks when there were none before, in the same way that $\bdd$-compressing a properly embedded annulus may create a compressing disk.  

For $X$ any compact manifold, let $|X|$ denote the number of components of $X$.  For example, if $S$ and $F$ are transverse compact surfaces, $|S \cap F|$ is the total number of arcs and circles in which $F$ and $S$ intersect.

\section{Bridges and bridge surfaces}

A properly embedded collection of arcs $T = \cup_{i=1}^n \alpha_i$ in a compact $3$-manifold is called {\em boundary parallel} if there is an embedded collection $E = \cup_{i=1}^n E_i$ of disks, so that, for each $1 \leq i \leq n$, $\bdd E_i$ is the end-point union of $\alpha_i$ and an arc in the boundary of the $3$-manifold.   If the manifold is a handlebody $A$, the arcs are called {\em bridges} and disks of parallelism are called {\em bridge disks}.

\begin{lemma}\label{lem:bddcompressible}
 Let $A$ be a handlebody and let $(T, \bdd T) \subset (A, \bdd A)$ be a collection of bridges in $A$. Suppose $F$ is a properly embedded surface in $A$ transverse to $T$ that is not a union of (punctured) disks and twice-punctured spheres.  If $F_T$ is incompressible in $A_T$ then $\bdd F  \neq \emptyset$ and $F_T$ is $\bdd$-compressible.

 \end{lemma}

 \begin{proof}  Suppose $D$ is a (punctured) disk component of $F$ or a twice-punctured sphere and let $F' = F - D$.  Suppose there is a compressing disk (resp. $\bdd$-compressing disk) for $F'_T$.  Then a standard cut and paste argument shows that there is a compressing disk (resp. $\bdd$-compressing disk) for $F'_T$ that is disjoint from $D$.  So with no loss of generality we may henceforth assume that {\em no} component of $F$ is a (punctured) disk or twice-punctured sphere.
 
A standard cut and paste argument provides a complete collection of meridian disks for $A$ which is disjoint from a complete collection of bridge disks for $T$.  Let the family $\Delta$ of disks be the union of the two collections, so in particular $A - \eta(\Delta)$ is a collection of $3$-balls.  In fact, choose such a collection, transverse to $F_T$, so that $|F_T \cap \Delta|$ is minimal.  $\Delta$ can't be disjoint from $F_T$, for then $F_T$ would be an incompressible surface in one of the ball components of $A - \eta(\Delta)$ and so a collection of disks, contrary to our assumption. If any component of  $F_T \cap \Delta$ were a closed curve, an innermost one on $\Delta$ would be inessential in $F_T$, since $F_T$ is incompressible, and an innermost inessential curve of intersection in $F_T$ (and perhaps more curves of intersection) could be eliminated by rechoosing $\Delta$.  We conclude that each component of $F_T \cap \Delta$ is an arc.  

Each arc in $F_T \cap \Delta$ can have ends on $\bdd F$, ends on $T$, or one end on both.  A similar cut and paste argument shows that if any arc has both ends on $\bdd F$, then an outermost such arc in $\Delta$ would be essential in $F_T$ and so the disk it cuts off in $\Delta$ would be a $\bdd$-compressing disk for $F_T$, as required.  Suppose then that all arcs of intersection have at least one end on $T$.  If any arc had both ends on $T$, a regular neighborhood of a disk cut off by an outermost such arc in $\Delta$ would contain a compressing disk for $F_T$ in its boundary, contradicting the hypothesis.  If all arcs of intersection have one end on $T$ and the other on $\bdd F_T$, a regular neighborhood of a disk cut off by an outermost such arc in $\Delta$ would contain a $\bdd$-compressing disk for $F_T$ in its boundary, as required.
  \end{proof}

Suppose $M$ is an irreducible $3$-manifold and $T$ is a properly embedded collection of arcs in $M$ transverse to a $2$-sphere $S$. Since $S$ bounds a $3$-ball it is separating, so $|S \cap T|$ is even.  In particular, if $|S \cap T| = 2$ then the core of the annulus $S_T$ can't bound a disk in $M_T$, for if it did the result would be a sphere in $M$ intersecting $T$ in a single point.  When $A$ is a handlebody and $T \subset A$ is a collection of bridges, more can be said:

\begin{lemma} \label{lem:irreducible}
    Suppose $A$ is a handlebody, $T$ is a collection of bridges
       in $A$ and $S \subset T $ is a sphere that intersects the link exactly
       twice transversally.  Then $T$ intersects the ball bounded by $S$ in an unknotted arc (ie a bridge in the ball). \end{lemma}

 \begin{proof}
We have seen that a $2$-sphere in $A$ intersects each bridge in an even number of points, so in particular $S$ intersects exactly one bridge $\alpha$.   Let $E$ be the bridge disk for $\alpha$. Since the core of the annulus $S_T$  can bound no disk in $M_T$ a standard
innermost disk argument allows $E$ to be chosen so that $S \cap E$ contains no closed curves at all.  Any arc of intersection between $S$
and $E$ must have both of its endpoints on $\alpha$ as $S \cap
\bdd A=\emptyset$. Thus there must be exactly one such arc of
intersection $\bbb$, cutting off a subarc $\alpha'$ of $\aaa$ on the opposite side of $S$ from $\bdd \aaa$.  Hence $\alpha'$ lies in the ball that $S$ bounds in $A$.   The subdisk of $E$ cut off by $\bbb$ is a bridge disk for $\alpha'$ in that ball.   \end{proof}

\begin{cor} Suppose $A$ is a handlebody, $T$ is a collection of bridges
       in $A$ and $D$ is a properly embedded (punctured) disk in $A$ whose boundary is inessential in $A_K$.  Then $D_K$ is isotopic rel boundary to a (punctured) disk $D'_K$ in $\bdd A$.  $D_K$ is punctured if and only if $D'_K$ is punctured.
       \end{cor}
       
       \begin{proof}  Apply Lemma \ref{lem:irreducible} to the sphere $D \cup_{\bdd D} D'$.  \end{proof}

Suppose $M$ is a closed $3$-manifold containing a link $K$.  A Heegaard surface $P$ for $M$ is a {\em bridge surface} for $K$ in $M$ if $K$ intersects each of the two complementary handlebodies of $P$ in $M$ in a collection of bridges.  

In trying to understand how two possibly different bridge surfaces are related, Lemma \ref{lem:irreducible} has the following pleasant corollary:

\begin{cor}  \label{cor:reminess} Suppose bridge surfaces $P$ and $Q$ intersect in a collection of curves such that any curve that is inessential in $P$ is also inessential in $Q$.  Then $Q$ can be properly isotoped, without adding curves of intersection and without removing any curve of intersection that is essential in both surfaces, until all curves of intersection are essential in $P$.
\end{cor}

\begin{proof}  The proof is by induction on the number of curves of intersection that are inessential in $P$.  If there are none, there is nothing to prove.

Among curves of intersection that are inessential in $P$, let $\gamma$ be innermost in $P$.  By hypothesis, $\gamma$ is also innessential in $Q$; let $D^P$ and $D^Q$ be the (punctured) disks that $\gamma$ bounds in, respectively $P$ and $Q$.  Then $S = D^P \cup_{\gamma} D^Q$ can be pushed slightly to be a sphere in one of the complementary handlebodies, say $X$, of $Q$.  Since $X$ is irreducible, the sphere $S$ bounds a ball $\mcB$ in $X$.  Since $K$ intersects each of $D^P$ and $D^Q$ in at most one point, it follows that $K$ is either disjoint from $\mcB$ or, following Lemma \ref{lem:irreducible}, intersects $\mcB$ in a single bridge.  In the latter case, a bridge disk for $\mcB$ can be isotoped so that it intersects $\gamma$ in a single point.    Thus in either case, $\mcB$ defines a proper isotopy of $D^Q$ to $D^P$ rel $\gamma$. Pushing $D^Q$ a  bit beyond $D^P$ removes $\gamma$ as a curve of intersection.  This completes the inductive step, once we show that any other curve of intersection that is removed by the isotopy is inessential in $Q$.  But since the interior of $\mcB$ is entirely disjoint from $Q$, the only other curves of intersection removed are those in $D^Q \cap P$.  Since these lie in $D^Q$ they are inessential in $Q$.  \end{proof}

Suppose both $P$ and $Q$ are bridge surfaces for a link $K$, so $M=A \cup_P B =X \cup_Q Y$ and the link $K$ is in bridge position with respect to both $P$ and $Q$.   That is, $K \cap A$, $K \cap B$, $K \cap X$ and $K \cap Y$ are all collections of bridges in the respective handlebodies.  

\begin{lemma} \label{lem:disjoint1}  Suppose $Q \subset A$ and $P_K$ compresses in $A_K - Q$.  Then either $K$ is the unknot in $S^3$ or $P_K$ is strongly compressible.
\end{lemma}

{\bf Remark:} In fact, we will show later (Lemma \ref{lem:unknot}) that the second alternative holds unless $K$ is in $1$-bridge position with respect to $P \cong S^2$.

\begin{proof}  With no loss of generality, suppose $P \subset X$, so $A \cap X$ is a cobordism between $P$ and $Q$ and the compressing disk for $P_K$ in $A_K - Q$ lies in $(A \cap X)_K$.  (See Figure \ref{fig:disjoint1a}.)  Let $P'$ be the surface obtained from $P$ by maximally compressing $P_K$ in $(A  \cap X)_K$.  $P'$ is a closed surface in the handlebody $X$ so it follows from Lemma \ref{lem:bddcompressible} that either $P'_K$ compresses in $X_K$ or $P'$ is the union of twice-punctured spheres. Moreover, by construction, $P'$ separates $Q$ from $P$.

    \begin{figure}[tbh]
    \centering
    \includegraphics[scale=0.5]{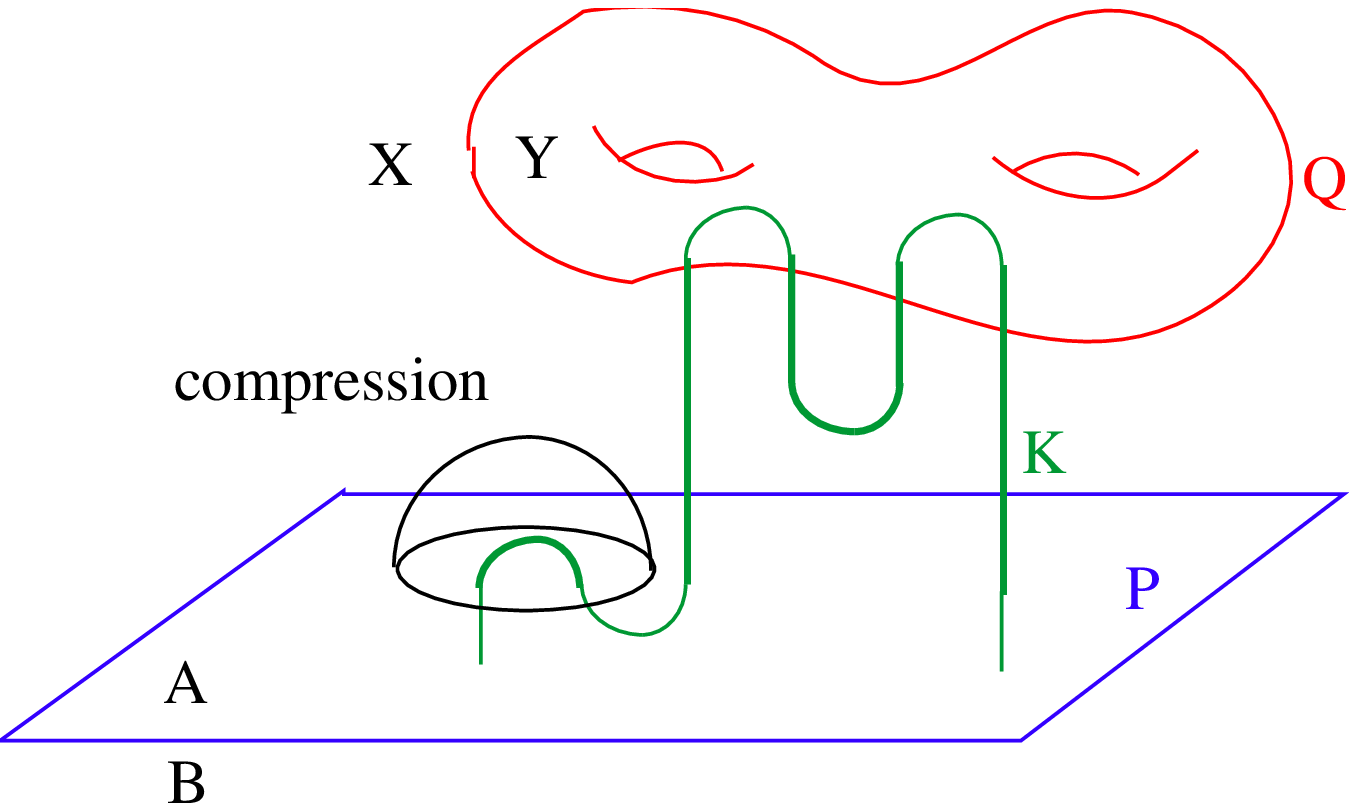}
    \caption{} \label{fig:disjoint1a}
    \end{figure}

Suppose first that $P'$ is a union of twice-punctured spheres.  Since $P'$ separates $Q$ from $P$, any arc from $P$ to $Q$ intersects $P'$ in an odd number of points.  Such an arc with a minimal number of intersection points with $P'$ will intersect each component of $P'$ in at most one point, so some component $P'_0$ of $P'$ (say the last component which the arc intersects before it intersects $Q$) will bound a ball in $A$ containing $Q$.  Dually, on the other side $P'_0$ bounds a ball in $X$ containing $P$.  By Lemma \ref{lem:irreducible} $P'_0$ is a bridge sphere for $K$ in $M$ so $K$ is the unknot in $S^3$.

    \begin{figure}[tbh]
    \centering
    \includegraphics[scale=0.5]{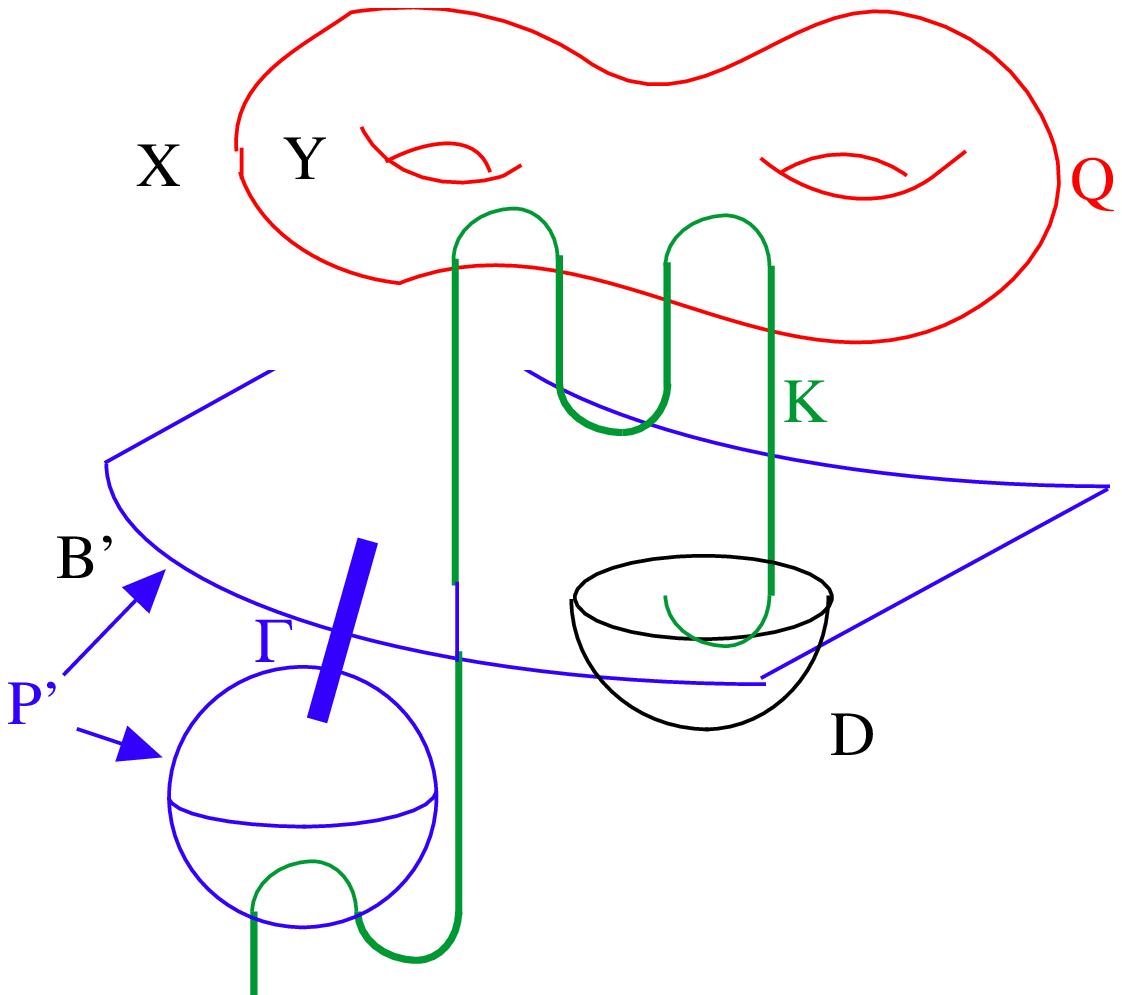}
    \caption{} \label{fig:disjoint1b}
    \end{figure}

Suppose on the other hand that $P'_K$ compresses in $X_K$.  By construction, any such compressing disk must lie on the side of $P'_K$ opposite to $Q$.  Denote that side $B'$ because it contains $B$.  In fact $B'$ is obtained from $B$ by attaching the $2$-handles determined by the maximal compression of $P_K$ in $(A  \cap X)_K$. Dually, $B$ is obtained from $B'$ by deleting a neighborhood of $P' \cup \Ggg$, where the graph $\Ggg$ consists of the arcs which are cocores of these $2$-handles and by construction $P' = \bdd B'$. Choose $\Ggg$ (up to slides and isotopies in $P' \cup \Ggg$) and choose the compressing disk $D \subset B'_K$ so as to minimize the number of points $|\Ggg \cap D|$. (See Figure \ref{fig:disjoint1b}.)  If $D$ is disjoint from $\Ggg$ then $D$ lies in $B_K$ and is disjoint from the compressing disks in $A_K$ dual to the arcs of $\Ggg$, exhibiting a strong compression of $P_K$.  

So we henceforth assume that $\Ggg$ intersects $D$.  Let $\Delta \subset B$ be a complete collection of meridian disks and bridge disks, chosen to minimize $|\Delta \cap D|$.  If $\Delta$ is a single bridge disk then $K$ is the unknot in $S^3$, as required.  So assume henceforth that $\Delta$ is more complicated.  In particular, each disk in $\Delta$ either is a compressing disk for $B_K$ or the boundary of its regular neighborhood is.  

Following \cite[Cor 2.3]{ST}, whose proof we now briefly recapitulate, consider the graph $\Upsilon$ in $D$ whose vertices are the points $\Ggg \cap D$ and whose edges are the arc components of $\Delta \cap D$.  As in that proof, the minimization of $|\Delta \cup D|$ guarantees that no closed component of $\Delta \cap D$ can bound a disk in the complement of $\Ggg$ and no loop in $\Upsilon$ cuts off a disk in $D$ that disjoint from $\Ggg$.  It follows that there is a vertex $v$ of $\Upsilon$ which is incident only to simple edges, perhaps because it is isolated.  If $v$ is not isolated then, of all the arcs incident to $v$, one that is outermost in $\Delta$ describes a way to slide the edge of $\Ggg$ containing $v$ to remove $v$ from $\Ggg \cap D$.  This would violate our original minimization of $|\Ggg \cap D|$. Hence $v$ is isolated, so $\Delta$ is disjoint from the compressing disk in $A_K$ dual to the edge containing $v$, exhibiting again a strong compression of $P_K$.
\end{proof}

\begin{lemma} \label{lem:disjoint2}  Suppose $Q \subset A$ and $P_K$ c-compresses in $A_K - Q$.  Then either $K$ is the unknot in $S^3$ or $P_K$ is c-strongly compressible.
\end{lemma}  

 \begin{proof}  

Let $D$ be a c-disk  for $A_K - Q$ in $P_K$.  If $D$ is disjoint from $K$, the result follows from Lemma \ref{lem:disjoint1}, so we may as well assume $D$ is a once-punctured disk.  As there, assume with no loss of generality $P \subset X$, so also $D \subset X$.  Let $\alpha \subset X$ be the bridge for $Q$ which intersects $D$ and let $E$ be a bridge disk for $\alpha$, so $\bdd E$ is the end-point union of $\alpha$ and an arc $\delta$ in $Q$.  Choose $E$ to minimize $|D \cap E|$.  If $D \cap E$ contains a closed curve then an innermost  such curve $c$ in $D$ cannot bound a once-punctured disk in $D$, else the union of this disk together with the disk that $c$ bounds in $E$ would be a once-punctured sphere in $X$, which is impossible.   Nor can $c$ bound a disk in $D_K$ else $c$ could be removed by a different choice of $E$.  Any arc component of $D \cap E$ with both end points on $\bdd D$ can be removed by rechoosing $E$.  We conclude that $D \cap E$ consists of a single arc, with one end at the point $K \cap D$ and the other end a point $p \in P \cap \bdd D$.  

    \begin{figure}[tbh]
    \centering
    \includegraphics[scale=0.5]{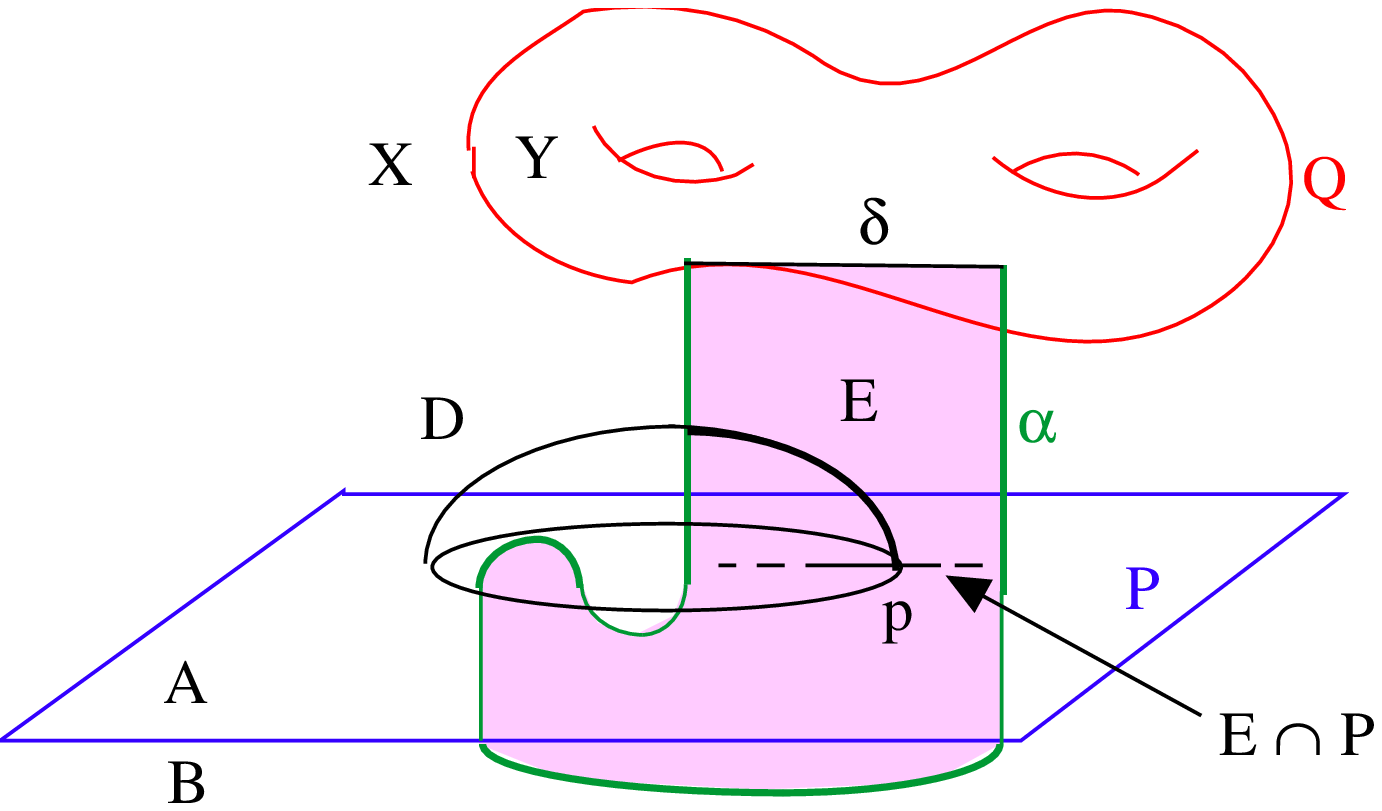}
    \caption{} \label{fig:disjoint2a}
    \end{figure}

Further choose $E$, retaining the requirement that $E \cap D$ is a single arc, to minimize $|P \cap E|$.  Notice that among the entire collection of curves $P \cap E$ only one curve intersects $\bdd D$, and it intersects it in the single point $\{ p \}$, since $\bdd D \cap E = \{ p \}.$  (See Figure \ref{fig:disjoint2a}.)  One possibility is that there is a closed curve $c$ of intersection that bounds a disk in $P_K$.  (Such a curve would necessarily intersect $\bdd D$ in an even number of points so, in particular, $c$ doesn't contain $p$.)  An innermost such curve in $P_K$ could be eliminated by a rechoice of $E$.  It follows that no closed curve of intersection bounds a disk in $P_K$.  If a closed curve of intersection $c \subset E \cap P$ that is innermost  in $E$ bounds a disk in $A_K$ then the result follows from Lemma \ref{lem:disjoint1}.  If $c$ bounds a disk $D'$ in $B_K$ and does not contain $p$ then $D$ and $D'$ are the disjoint disks that are sought. If $c$ does contain $p$ then the boundaries of $D$ and $D'$ intersect precisely in the single point $p$.  The boundary of a regular neighborhood of $\bdd D \cup \bdd D'$ in $P$ is essential in $P_K$ (since $|K \cap P| \geq 2$) and also bounds a disk $D''$ in $B_K$, namely two copies of $D'$ banded together along $\bdd D$.  Then $D$ (possibly punctured) and $D''$ are the disjoint disks that are sought.  So henceforth we assume that there are no closed curves in $P \cap E$.

    \begin{figure}[tbh]
    \centering
    \includegraphics[scale=0.5]{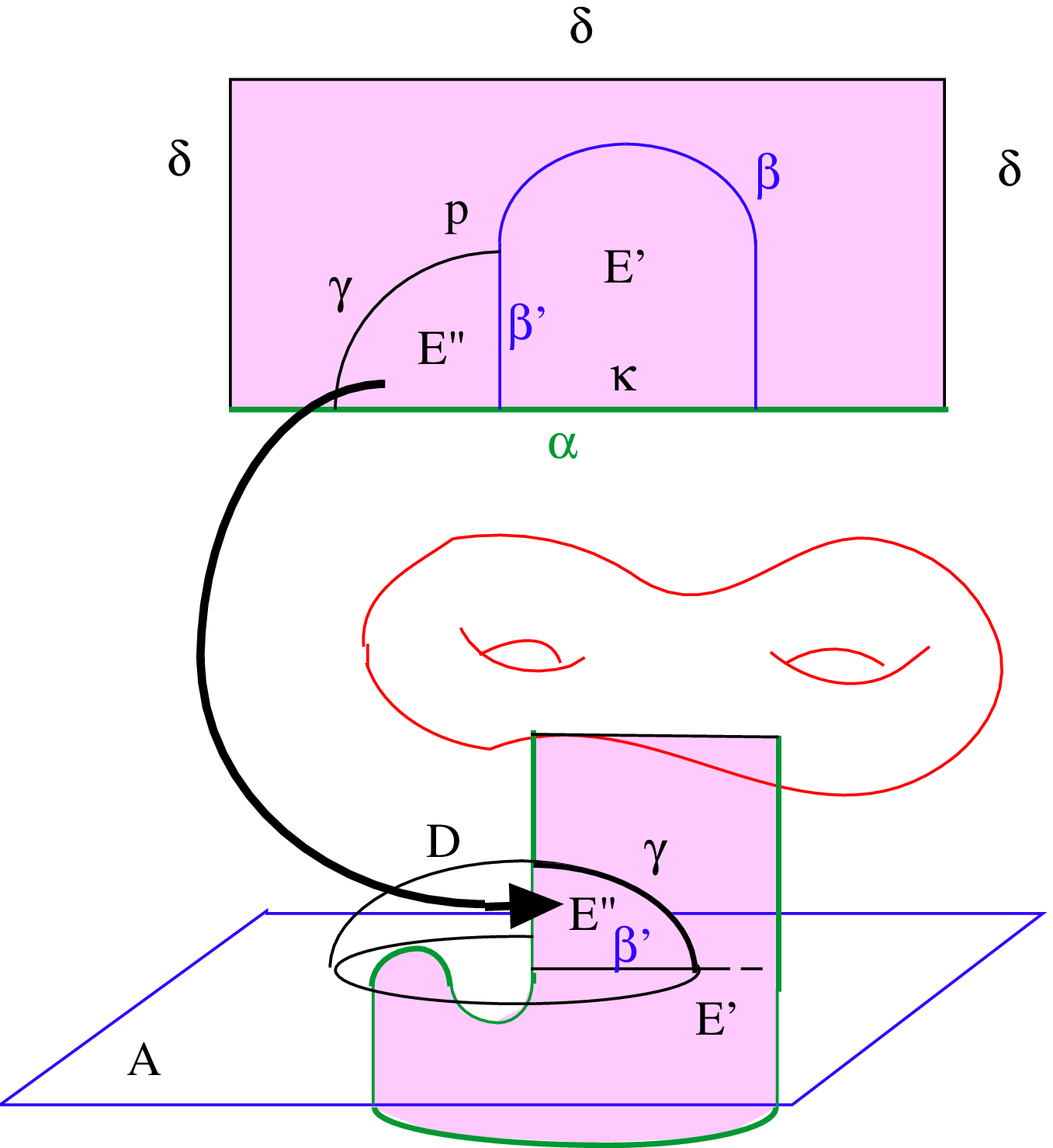}
    \caption{} \label{fig:disjoint2b}
    \end{figure}

Consider then the arcs of intersection in $P \cap E$.  $\bdd E$ is the end-point union of two arcs, $\aaa$ a component of $K - Q$ and $\ddd$ an arc in $Q$.  Since $Q \subset A$, $\bdd E$ intersects $P$ only in points on $\aaa$.  An outermost arc of intersection $\bbb$ in $E$ then cuts off a bridge disk $E'$ for a bridge $\kkk \subset \aaa$ for $K$ with respect to $P$.  In particular, a regular neighborhood of $E'$ contains in its boundary a compressing disk $D'$ for $P_K$.  If $D' \subset A$ we are done by Lemma  \ref{lem:disjoint1}.  If $D' \subset B$ and $\bbb$ does not contain $p$,  $D'$ and $D$ are the required disjoint c-disks for $P_K$.  The only remaining possibility is that $p \in \bbb$, $\bbb$ is the only outermost arc of intersection and the disk $E'$ it cuts off lies in $B_K$.  In this case, consider the arc $\gggg = D \cap E \subset A_K$ with one end on $\aaa$ and the other end on $p$.  Then the union of $\gggg$ and a subarc $\bbb'$ of $\bbb$ cut off a subdisk $E'' \subset A$ of $E - E'$.  (See Figure \ref{fig:disjoint2b}.)  The union of $D$ and $E''$ along $\gggg$ has a regular neighborhood in $A$ consisting of two disks, one parallel to $D$ in $A_K$ and the other, $D''$, parallel to $D$ in $A$ but disjoint from $K$.  $\bdd D''$ must be essential in $P_K$, for if it were inessential then it would bound a disk in $P_K$ and so $\bdd D$ would bound a once-punctured disk in $P_K$, contradicting the assumption that $D$ is a c-disk.  Hence $D''$ is a compressing disk for $P_K$ and the result once again follows from Lemma   \ref{lem:disjoint1}. 
\end{proof}

The hypotheses can be weakened further, to allow inessential curves in $P_K \cap Q_K$:

\begin{lemma} \label{lem:disjoint3}  Suppose 
\begin{itemize}
\item every curve in $P_K \cap Q_K$ is inessential in both $P_K$ and $Q_K$
\item $Q_K \cap A$ contains some curve that is essential in $Q_K$ and 
\item $P_K$ c-compresses in $A_K - Q$.
\end{itemize}  Then either $K$ is the unknot in $S^3$ or $P_K$ is c-strongly compressible.   
\end{lemma}  

\begin{proof}   The proof is by induction on $|P_K \cap Q_K|$.  If $P_K$ and $Q_K$ are disjoint then the result follows from Lemma \ref{lem:disjoint2}.  (Note that if $Q_K$ contains no essential curves then $Q_K$ is a twice-punctured sphere, so $K$ is the unknot in $S^3$.)  Suppose, for the inductive step, that  $|P_K \cap Q_K| \geq 1$ and $c$ is a circle of intersection that is innermost in $P$.  Let $D \subset P$ be the (punctured) disk that $c$ bounds in $P$ and $E$ be the (punctured) disk that $c$ bounds in $Q$.   Clearly the c-disk for $P_K$ disjoint from $Q_K$ is disjoint from $D$.  Although the interior of $E$ may intersect $P$, the interior of $D$ is disjoint from $Q$, so after a slight push on $E$,  $S = D \cap_c E$ is a sphere that is disjoint from $Q$.  Hence $S$ bounds a ball in $X$ or $Y$.  $S$ itself is either disjoint from $K$ or punctured twice.  So $S$ either bounds a ball in $X_K$ or, according to Lemma \ref{lem:irreducible}, it bounds a ball in $X$ which $K$ intersects in an unknotted arc.  In either case, the ball describes a proper isotopy of $Q_K$, stationary away from $E$, that replaces $E$ with $D$.  After a further small push, $c$ is removed as a curve of intersection.  That is, $|P_K \cap Q_K|$ has been reduced by at least one.  Any essential curve in $Q_K$ that is disjoint from $P_K$ is clearly unaffected by this isotopy, since the curve must be disjoint from $E$, so the curve remains in $A$. Furthermore, the $P_K$ compressing disk in $A_K$ disjoint from $Q$ remains disjoint from $Q$, since it was disjoint from $D$.   This completes the inductive step.  \end{proof}

\section{Conway spheres and bridge position} \label{sec:Conway}

\begin{defin}  \label{def:Conwaysphere}  A {\em Conway sphere} (\cite{Co}) $S$ for a link $K$ in a $3$-manifold $M$ is a sphere $S \subset M$ transverse to $K$ such that $|K \cap S| = 4$.  A Conway sphere is an {\em incompressible Conway sphere} if $S_K$ is incompressible in $M_K$. 
\end{defin}

   \begin{figure}[tbh]
    \centering
    \includegraphics[scale=0.5]{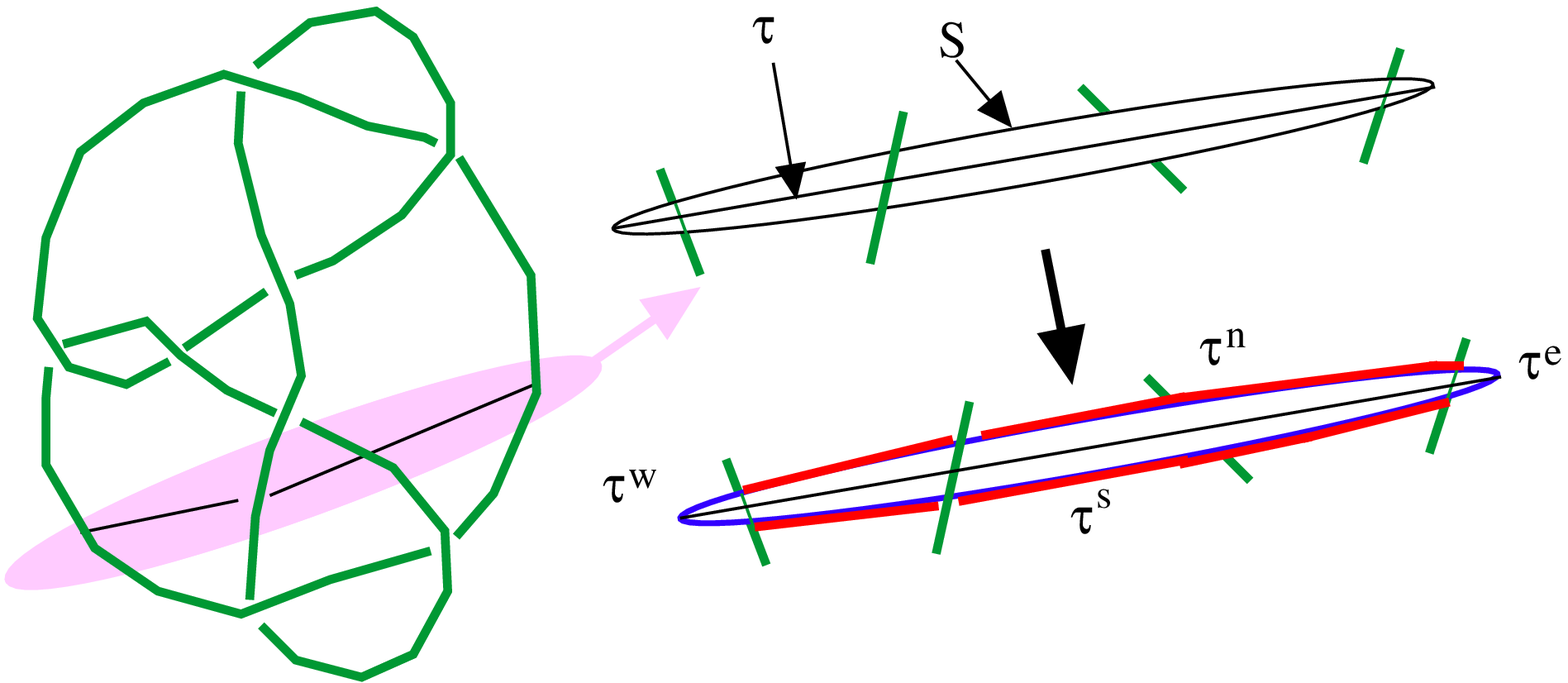}
    \caption{} \label{fig:Conwaysphere}
    \end{figure}

Suppose $M_0, M_1$ are orientable $3$-manifolds containing links $K_0$, $K_1$ respectively.  For each $i = 0, 1$, let $\tau_i$ be an arc in $M_i$ whose ends lie on $K_i$ but $\tau_i$ is otherwise is disjoint from  $K_i$.  Let $\mcB_i$ be a regular neighborhood of $\tau_i$, a ball intersecting $K_i$ in two arcs, one near either end of $\tau_i$.  Then $S_i = \bdd \mcB_i$ is a Conway sphere for $K_i$ in $M_i$.   The arcs $K_i \cap B_i$ are parallel in $\mcB_i$ to arcs $\tau^e_i$ and $\tau^w_i$ in $S_i$.  Let $\tau^n_i$ and $\tau^s_i$ denote a pair of arcs in each $S_i$ which, together with $\tau^e_i$ and $\tau^w_i$, form an embedded circle in $S_i$.  See Figure \ref{fig:Conwaysphere}.  

\begin{defin}  \label{def:tangleop} Given $\tau_i$ as above, let $K_0 +_c K_1$, called a {\em Conway sum} of the $K_i$, denote a link in $M_0 \# M_1$ obtained by removing the interior of $\mcB_i$ from each $M_i$ and gluing $S_0$ to $S_1$ via a homeomorphism that identifies the pair of arcs $\tau^e_0, \tau^w_0$ with the arcs $\tau^e_1, \tau^w_1$.  

Similarly, let $K_0  \times_c K_1$, called a {\em Conway product}, denote a link in $M_0 \# M_1$ obtained by instead gluing $S_0$ to $S_1$ via a homeomorphism that identifies the pair of arcs $\tau^e_0, \tau^w_0$ with the arcs $\tau^n_1, \tau^s_1$ and the pair of arcs $\tau^n_0, \tau^s_0$ with the arcs $\tau^e_1, \tau^w_1$. 

The image $S$ of $S_0$ and $S_1$ after their identification is called the Conway sphere of the sum (or product).
\end{defin}

\bigskip

{\bf Note:}  Essentially the same constructions can be done for disjoint arcs $\tau_0, \tau_1$ that are contained in the same manifold $M$ and which have their ends on the same link $K \subset M$.  In that case choose the identification $S_0 \cong S_1$ to be orientation reversing, so the resulting manifold is $M \# S^1 \times S^2$.  The sum and product links in $M \# S^1 \times S^2$ respectively are denoted $K +_c$ and $K  \times_c$ and in this case the Conway sphere of the sum or product is non-separating.  

\bigskip

The definition of Conway sum and product is motivated by analogy to tangle sum and tangle product (cf \cite[p 47-48]{Ad}).  Unlike the standard construction of a connected-sum of links, a Conway sum or Conway product of links depends on many choices beyond the question of which components of $K_i$ contain the ends of $\tau_i$. Most prominent is the choice of the arcs $\tau_i$ defining the operation, but there is also some choice in how $S_0$ is identified to $S_1$ beyond the constraints given by the definitions.  

In general, one would expect little connection between bridge presentations of $K_i$ in $M_i$ and bridge presentations of $K_0 +_c K_1$ and $K_0 \times_c K_1$ in $M_0 \# M_1$ (or $K \times_c$ in $M \# S^1 \times S^2$).  The most obvious problem is that each $\tau_i$ may intersect a bridge surface $P_i$ for $K_i$ an unknown number of times, eg possibly $|\tau_0 \cap P_0| \neq |\tau_1 \cap P_1|$, and there is no way of cobbling together the complementary punctured bridge surfaces into a plausible bridge surface for the resulting link.  In the case of Conway sum, the problem can be alleviated by limiting the operation to arcs $\tau_i$ that lie in the bridge surface $P_i$, and requiring that the equator curve $S_0 \cap P_0$ be identified with the equator curve $S_1 \cap P_1$.  Then the resulting surface $P_0 \# P_1$ is the standard connected sum of Heegaard surfaces for $M_0$ and $M_1$, so $P_0 \# P_1$ is a Heegaard surface for $M_0 \# M_1$ and so a potential bridge surface for $K_0 +_c K_1$.  But even in this case, if the arcs $\tau_i$ intersect bridge disks for $K_i$, then there is no natural reason why $K_0 +_c K_1$ should be in bridge position with respect to $P_0 \# P_1$.  So one expects that Conway sums do not in general behave well with respect to bridge number.  

\bigskip

The situation is more hopeful for Conway products:

\begin{defin}  \label{def:conwayprod} Suppose, for $i = 0, 1$, $P_i$ is a bridge surface for the link $K_i$ in $M_i$ and $\tau_i$ is an arc in $P_i$ intersecting $K_i$ exactly in the end points of $\tau_i$.  Form a Conway product by taking $\tau^n_i, \tau^s_i$ to be arcs in $S_i$ disjoint from  the equator $\bdd P_i$.  The result is called a Conway product that {\em respects the bridge surfaces}.  The sphere $S$ is called a {\em Conway decomposing sphere} for the pair $(P_0 \# P_1, K_0 \times_c K_1)$.

The same terminology is used when $P$ is a bridge surface for the link $K$ in $M$;  $\tau_0$ and $\tau_1$ are disjoint arcs in $P$ intersecting $K$ exactly in their end points; and a Conway product is formed by taking $\tau^n_i, \tau^s_i$ to be arcs in $S_i$ disjoint from  the equator $\bdd P_i$.  $S$ is then a Conway decomposing sphere for the pair $(P \# (S^1 \times S^1), K \times_c)$.
\end{defin}

The reason for regarding $S$ as a decomposing sphere for the pair is this:

\begin{prop} \label{prop:conwayprod} Let $K_0 \times_c K_1$ be a Conway product respecting the bridge surfaces $P_i$ for $K_i$ in $M_i, i = 0, 1$.  Then $P_0 \# P_1$ is a bridge surface for $K_0 \times_c K_1$ in $M_0 \# M_1$.  Moreover the bridge numbers satisfy $$\beta (K_0 \times_c K_1) \leq \beta (K_0) + \beta (K_1) - 1.$$

Similarly, when $K \times_c$ is a Conway product respecting the bridge surface $P$ for $K$ in $M$ then $P \# (S^1 \times S^1)$ is a bridge surface for $K \times_c$ in $M \# (S^1 \times S^2)$ and the bridge number satisfies $\beta (K \times_c) \leq \beta (K) - 1.$

\end{prop}

\begin{proof}  We consider only the first case, when the Conway product $K = K_0 \times_c K_1$ is of links in two different manifolds; the proof for a product $K \times_c$ is similar.

Examine how the bridge surfaces $P_0$ and $P_1$ intersect $S_0$, after $S_0$ and $S_1$ are identified.  Each $P_i$ intersects $S_0$ in a single closed curve $c_i = S_0 \cap P_i$.  Since $\tau_0$ lies in $P_0$, $c_0$ intersects each of $\tau^e_0, \tau^w_0$ in a single point and, by hypothesis, is disjoint from $\tau^n_0, \tau^s_0$.  Similarly $c_1$ intersects each of $\tau^n_0, \tau^s_0$ in a single point and is disjoint from each of $\tau^e_0, \tau^w_0$.  It follows that the $c_i$ can be isotoped in $S_0 = S_1 = S$ rel the circle $\overline{\tau} = \tau^e_0 \cup \tau^n_0 \cup \tau^w_0 \cup \tau^s_0$ so that $c_0 \cap c_1$ consists of a single point in each of the disks $S - \overline{\tau}$ and together the curves $c_0 \cup c_1$ divides $S$ into quadrants, each containing a single puncture.  The three curves $c_0 \cup c_1 \cup \overline{\tau}$ divide $S$ into octants.  See Figure \ref{fig:Conwayproda}.

  \begin{figure}[tbh]
    \centering
    \includegraphics[scale=0.5]{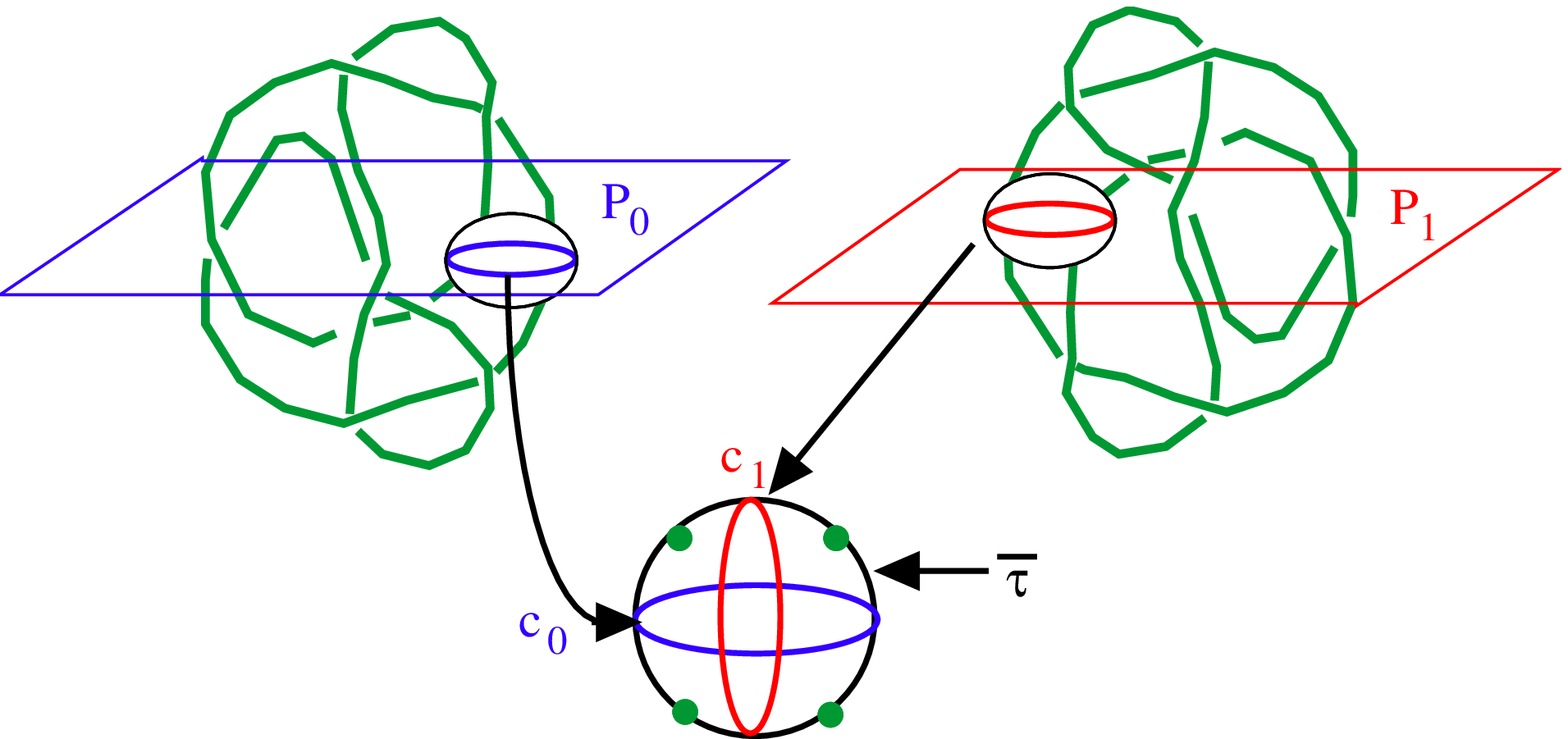}
    \caption{} \label{fig:Conwayproda}
    \end{figure}

With no loss of generality, assume $\overline{\tau}$ and the $c_i$ are mutually orthogonal great circles of $S$.  Let $\rho_t: S^1 \to S$ be an isotopy of $c_0$ that rotates $c_0$ through an angle of $\pi/2$ to $c_1$ around the two fixed points $c_0 \cap c_1$.  Thus $\rho_i (S^1) = c_i, i = 0, 1$.  During the isotopy the image of $\rho_t$ will cross exactly two of the four punctures $S \cap K$.  (There are two choices for such an isotopy, one for each direction of rotation; opposite choices will carry $c_0$ to anti-parallel copies of $c_1$ and will cross the complementary pair of punctures.)  

Let $S \times I$ be a collar of $S$ in $M_0 \# M_1$ and use it to consider an alternate construction of $M_0 \# M_1$, namely identify $S_i = \bdd M_i$ with $S \times \{ i \} \subset S \times I, i = 0, 1$.  Connect $\bdd (P_0 - \mcB_0) = c_0$ to $\bdd (P_1 - \mcB_1)= c_1$ in $S \times I$ via the embedded annulus $S^1 \times I \subset S \times I$ given by $(x, t) \to (\rho_t(x) , t)$ and connect the four points $\bdd (K_0- \mcB_0)$ to $\bdd (K_1- \mcB_1)$ in $S \times I$ by the $4$ arcs $(K \cap S) \times I$.  See Figure \ref{fig:Conwayprodb}. The result is a specific embedding of  the surface $P_0 \# P_1$ and $K$ in $M_0 \# M_1$ so that $|(P_0 \# P_1) \cap K| = |P_0  \cap K_0| + |P_1 \cap K_1| - 2$.  So it remains only to show that $P_0 \# P_1$ is a bridge surface for $K$ in $M_0 \# M_1$.  Part of this is easy: as observed earlier, $P_0 \# P_1$ is the standard connected sum of the Heegaard surfaces $P_i$ for $M_i$, so it is a Heegaard surface $M_0 \# M_1$.  The remaining problem then is to exhibit bridge disks for each arc component of $K - (P_0 \# P_1).$  This requires a more concrete description:

 \begin{figure}[tbh]
    \centering
    \includegraphics[scale=0.5]{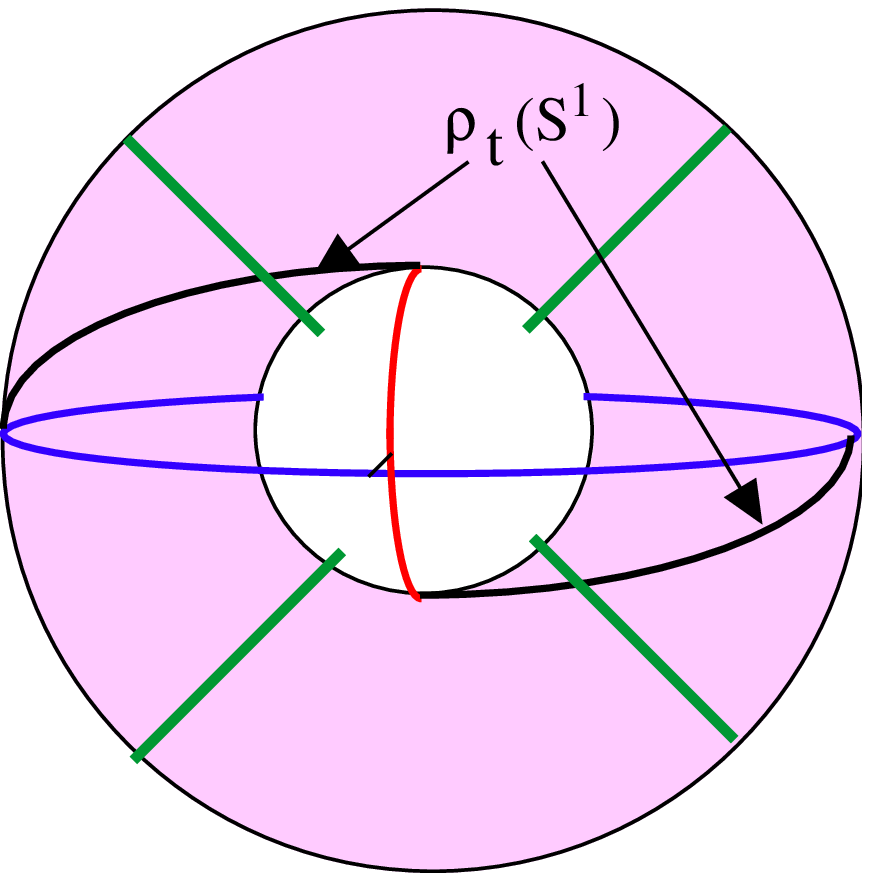}
    \caption{} \label{fig:Conwayprodb}
    \end{figure}

As is standard in the theory of tangles (cf \cite{GL}), label the four points of $K \cap S$ according to the four quadrants in which we imagine them lying: $NW, NE, SE, SW$.  With no loss of generality, assume $c_0$ is horizontal, separating the northern hemisphere from the southern, and so separating the pair $NW, NE$ from the pair $SE, SW$.  Similarly assume that $c_1$ is vertical, separating the western hemisphere from the eastern, and so separating the pair $NW, SW$ from the pair $NE, SE$.  Further assume that $\rho_t$ isotopes $c_0$ across the points $NW$ and $SE$ to $c_1$.   

In order to exhibit bridge disks, consider the types of bridges (ie arc components of $K - (P_0 \# P_1)$) that can arise.  First consider a bridge that is disjoint from $S \times I$, eg a component $\alpha$ of $K_0 - P_0$ disjoint from $S_0$.  Let $E$ be a bridge disk for $\alpha$ with respect to the original splitting surface $P_0$.  By general position, $\bdd E$ intersects the original arc $t_0 \subset P_0$ transversally in a number of points.  Corresponding to each point is a half-meridian disk (denoted $\mu_0$ in Figure \ref{fig:Conwayprodc} in which $E$ intersects $\mcB_0 - P_0$.  Thus $E$ intersects $S_0$ in a collection of parallel arcs in one of the twice-punctured hemisphere components (say the northern hemisphere) of $S_0 - P_0 \cong S - c_0$.  Push all these arcs $E \cap (S_0 - K)$ to the west of $c_1$ in the northern hemisphere.  The isotopy $\rho_t$ sweeps $c_0$ entirely across all the arcs $E \cap (S_0 - K)$ and, in doing so, defines a collection of disks which replace the half meridians of $\tau_0$ to give a bridge disk for $\alpha$.  

 \begin{figure}[tbh]
    \centering
    \includegraphics[scale=0.5]{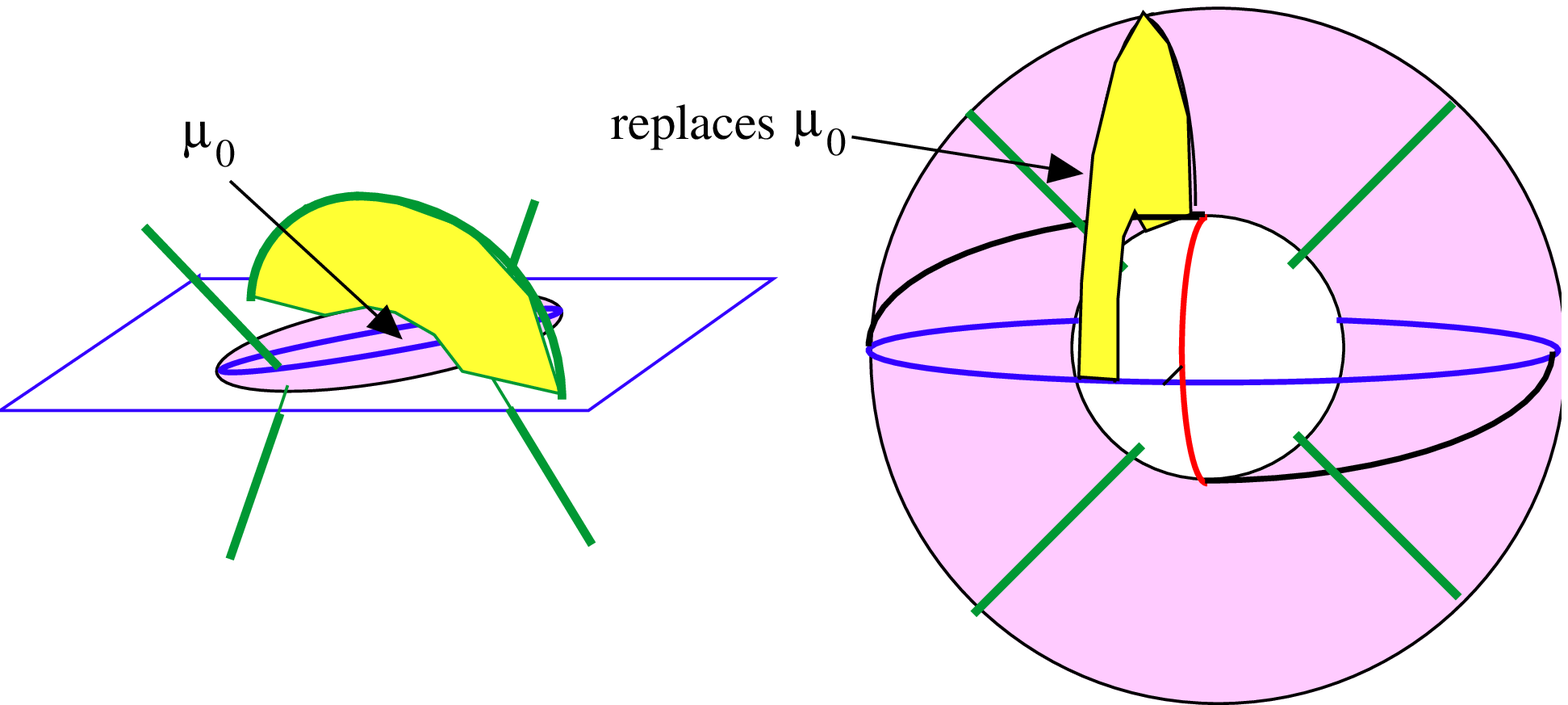}
    \caption{} \label{fig:Conwayprodc}
    \end{figure}

The argument is little different for a bridge $\alpha$ of $K - (P_0 \# P_1)$ that contains exactly one of the four points $NW \times \bdd I$ or $SE \times \bdd I$, say $NW \times \{ 0 \} \in S \times  \{ 0 \}$. Let $E$ be the bridge disk in $M_0 - P_0$ for the bridge $\alpha_0$ that contains the arc .   Then the northern hemisphere of $S_0 - P_0$ cuts off from $E$ a collection of half-meridians of $\tau_0$ just as before, together with a disk whose boundary consists of three arcs, one of them lying on $P_0 \cap \mcB_0$, one of them the segment $\alpha_0 \cap \mcB_0$, and one of them an arc $\gamma$ from $NW$ to $c_0$ in $S \times \{ 0 \}$ which we may as well take to be part of $\tau^w_0$.  But the sweep of $c_0$ across $\gamma$ defines a disk that can be used (together with the half-meridians of the previous case) to complete $E - \mcB_0$ to give a bridge disk for $\alpha$.  See Figure \ref{fig:Conwayprodd}.

 \begin{figure}[tbh]
    \centering
    \includegraphics[scale=0.5]{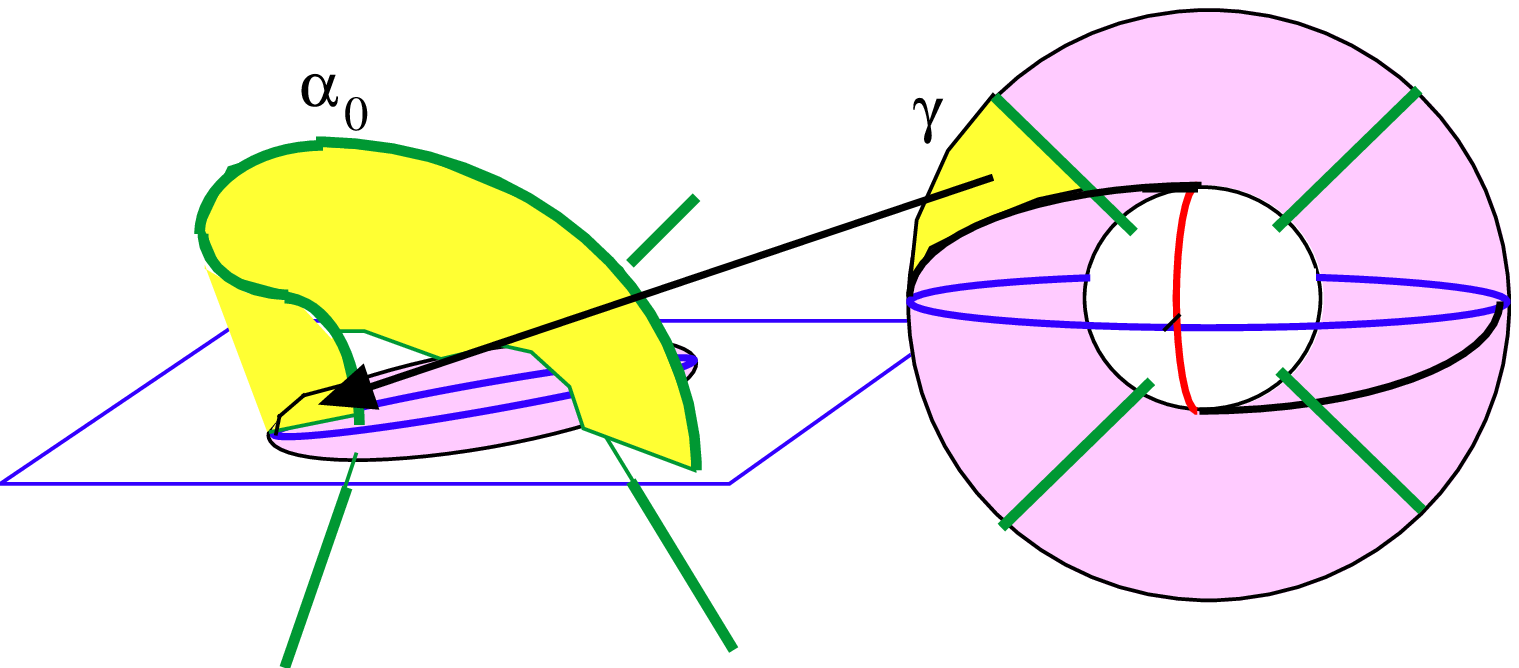}
    \caption{} \label{fig:Conwayprodd}
    \end{figure}

Next consider the case when $\alpha$ intersects $S \times \{ 0 \}$ in exactly one point, but a point that the isotopy $\rho_t$ does not sweep across, say the point $NE \times \{ 0 \}$.  Then all of the arc $NE \times I \subset S \times I$ lies in $\alpha$.  Assuming (for the moment) that $\alpha$ otherwise does not intersect $S \times \bdd I$, $\alpha$ is then the union of $NE \times I$ and two other segments: $\alpha_0$ and $\alpha_1$, where $\alpha_i = \alpha \cap (M_i - \mcB_i)$.  Let $E_i$ be a bridge disk for the bridge in $K_i - P_i$ that contains $\alpha_i$.  Much as before, we can arrange that $E_0$ intersects $S$ in a collection of arcs (with both ends on $c_0$) in the western half of the northern hemisphere, together with a single arc $\gamma$ connecting $NE$ to $c_0$.  It will be useful in this phase to take for $\gamma$ an arc whose other end is at one of the two rotation points $c_0 \cap c_1$.   See Figure \ref{fig:Conwayprode}.   Dually, $E_1$ intersects $S$ in a collection of arcs (with both ends on $c_1$) in the eastern half of the southern hemisphere, together with an arc connecting $NE$ to $c_1$, for example $\gamma$.  The arc $\gamma$ is unaffected by the isotopy of $c_0$ to $c_1$, so $\gamma \times I$ attaches $E_0 - \mcB_0$ to $E_1 - \mcB_1$, creating a single bridge disk for $\alpha$.

 \begin{figure}[tbh]
    \centering
    \includegraphics[scale=0.5]{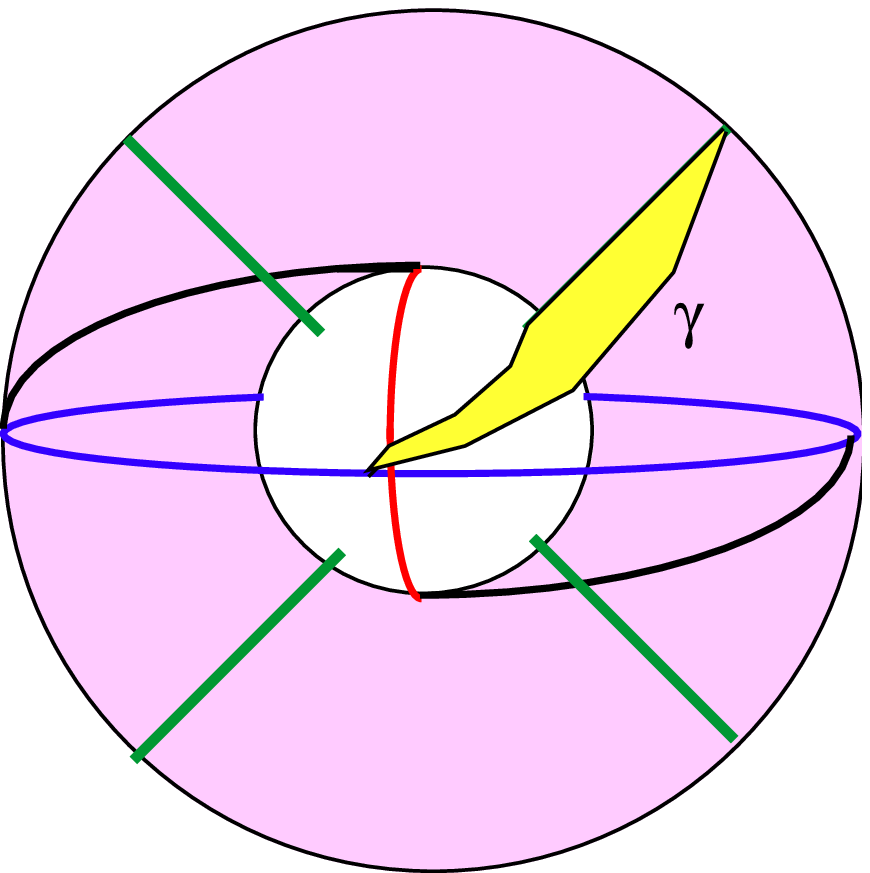}
    \caption{} \label{fig:Conwayprode}
    \end{figure}

The two other points of $S \times \bdd I$ that the bridge $\alpha$ containing $NE \times I$ could contain are the points $NW \times \{ 0 \}$ and $SE \times \{ 1 \}$.  In that case, the bridge disk for $\alpha$ is assembled by combining the arguments above.  
\end{proof}  

\begin{figure}[tbh]
    \centering
    \includegraphics[scale=0.5]{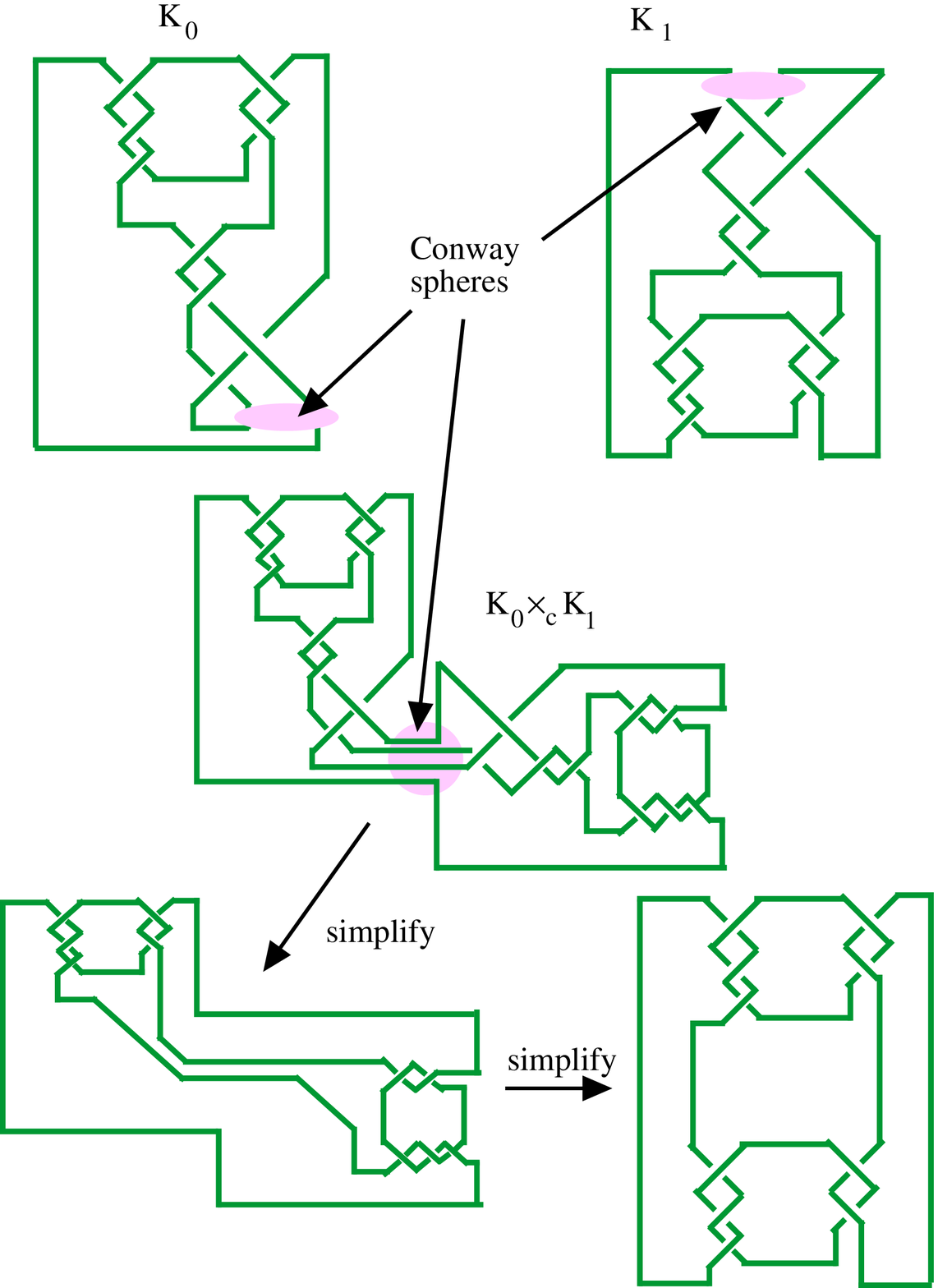}
    \caption{} \label{fig:blairfinal}
    \end{figure}

It's natural to ask whether the inequality in Proposition \ref{prop:conwayprod} is an equality, just as Schubert showed the analogous equality for connected sum.  (See \cite{Schub}, or \cite{Schul} for a modern proof.)  Ryan Blair pointed out an example in $S^3$ where it is not.  The factor knots are $3$-bridge knots that incorporate rational tangles that cancel when the product is constructed.  The resulting Conway product is a $4$-bridge link.  See Figure \ref{fig:blairfinal}.

\section{Spines and sweep-outs}

A spine of a handlebody $A$ is a properly embedded finite graph $\Sigma$ in $A$ (typically chosen to have no valence $1$ vertices) so that $A - \Sigma \cong \bdd A \times [0, 1)$.  Given a spine $\Sigma$ and a collection $T$ of bridges in $A$, $T$ can be isotoped in $A$ (for example by shrinking a collection $E$ of bridge disks very close to $\bdd A$) so that the projection (called the height) $A - \Sigma \cong \bdd A \times [0, 1) \to  [0, 1)$ has a single maximum on each bridge $\aaa_i$.  For each $\alpha_i$ connect $\Sigma$ to that maximum by an arc in $A$ which is monotonic with respect to height.    The union of $\Sigma$ with that collection of arcs is called a spine $\Sigma_{(A,T)}$ of $(A, T)$.  Note that there is a homeomorphism $A - \Sigma_{(A,T)} \cong \bdd A \times [0, 1)$ which carries $T - \Sigma_{(A,T)}$ to $(\bdd A \cap T) \times [0, 1)$.  Put another way, there is a map $(\bdd A, \bdd A \cap T) \times I \to (A, T)$ which is a homeomorphism except over $\Sigma_{(A,T)}$, and the map gives a neighborhood of $\Sigma_{(A,T)}$ a mapping cylinder structure.

Suppose that a link $K \subset M$ is in bridge position with respect to a Heegaard surface $P$ for $M$.  Then the closed complementary components of $P$ are handlebodies $A$ and $B$ that $K$ intersects in a collection of bridges.  Let $\Sigma_{(A,K)}$ (resp $\Sigma_{(B,K)}$) denote a spine in $A$ (resp $B$) for $K \cap A$ (resp $K \cap B$).  Then, following the above remarks,  there is a map $H:
(P, P \cap K) \times I \to (M, K)$ that is a homeomorphism except over
$\Sigma_{(A,K)}\cup \Sigma_{(B,K)}$ and, near $P \times \bdd I$, the map $H$ gives 
a mapping cylinder structure to a neighborhood of $\Sigma_{(A,K)} \cup \Sigma_{(B,K)}$. Little is lost and some brevity gained if we restrict $H$ to $P_K \times (I, \bdd I) \to (M_K, \Sigma_{(A,K)}\cup \Sigma_{(B,K)})$.  $H$ is then called a {\em sweep-out} associated to $P$.  

As a warm up, here is a classical application of sweep-outs.  In effect it strengthens somewhat the conclusions of Lemmas \ref{lem:disjoint1}, \ref{lem:disjoint2} and \ref{lem:disjoint3}, for it shows, in those lemmas, that if $P_K$ is c-weakly incompressible then $K$ is not only the unknot in $S^3$ but it is also in $1$-bridge position with respect to $P$.  

\begin{lemma} \label{lem:unknot}  Suppose $K \subset S^3$ is the unknot and is in bridge position with respect to $P$.  Either $P$ is a sphere and $K$ is the unknot, in $1$-bridge position with respect to $P$, or $P_K$ is strongly compressible.
\end{lemma}  

\begin{proof}   We will assume that either the number of bridges $m \geq 2$ or $P \neq S^2$, for otherwise $K$ is $1$-bridge with respect to $P$ a sphere and we are done.  In particular, we can assume that the neighborhood of any bridge disk  for $K$ contains a compressing disk for $P_K$ in $S^3_K$.   

Choose spines $\Sigma_{(A,K)}$ and $\Sigma_{(B,K)}$ as described above.  In particular, the end points of $\Sigma_{(A,K)}$ (resp $\Sigma_{(B,K)}$) are just the collection of $m$ maxima (resp minima) of $K$ with respect to a sweep-out $H: P_K \times (I, \bdd I) \to (S^3_K, \Sigma_{(A,K)}\cup \Sigma_{(B,K)})$.   Let $D$ be the disk that $K$ bounds.  For $\epsilon$ very small, $P^{\epsilon}_K = H(P_K \times \{ \epsilon \})$ is the boundary of a regular neighborhood of $\Sigma_{(B,K)}$.  By transversality (of $D$ with $\Sigma_{(B,K)}$), the curves of intersection $D \cap P^{\epsilon}_K$ consist of a family of $m$ arcs (the family is $\bdd$-parallel in $D$ and each corresponds to a minimum of $K$) and an unnested collection of simple closed curves, each corresponding to a point in $D \cap \Sigma_{(B,K)}$.   Each arc and simple closed curve cuts off a disk from $D$; the former are bridge disks for $\Sigma_{(B,K)}$ and the latter are compressing disks for $P_K$ in $B_K$.  

Similar comments hold for $P^{1 - \epsilon}_K = H(P_K \times \{ 1 - \epsilon \})$, except all the subdisks of $D$ cut off lie in $A_K$.  Now consider the intersection with $D$ of a generic $P^{t}_K = H(P_K \times \{ t \})$.  Any circle or arc in $D$ is clearly inessential in $D$ and $P^{t}_K$ necessarily intersects $D$ in $m$ arcs.  After all circles in $D \cap P^{t}_K$ that are inessential in $P^t_K$ are removed by an isotopy of $D$, either an innermost circle of intersection or an outermost arc of intersection cuts off a subdisk of $D$ that is either a bridge disk or a compressing disk for $P^t_K$.  Furthermore the neighborhood of a bridge disk contains a compressing disk for $P^t_K$.  So we conclude that for any generic $t$ there is some curve of intersection that defines a compressing disk for $P^t_K$,  either lying in $A_K$ or in $B_K$.  

Since for small $t$ such a compressing disk lies in $B^K$ and for large $t$ such a disk lies in $A_K$ and always there is some compressing disk, it follows that there is a (possibly non-generic) $t$ at which one has compressing disks for $P_K$ both in $A_K$ and $B_K$.  Since these two disks are defined by disjoint arcs of intersection, the disks themselves will be disjoint, and so define a strong compression of $P_K$.  (A short and standard argument shows this is true even when $t$ is non-generic, cf Lemma \ref{lem:tworegion} below.)
\end{proof}

Suppose both $P$ and $Q$ are possibly different bridge surfaces for the link, so $M=A \cup_P B =X \cup_Q Y$ and the link $K$ is in bridge position with respect to both $P$ and $Q$.  Choose spines $\Sigma_{(A,K)},  \Sigma_{(B,K)}, \Sigma_{(X,K)}, \Sigma_{(Y,K)}$ in each handlebody, all in general position (hence disjoint) from each other in $M$.  Using these spines, as above, choose sweepouts for $P_K$ and $Q_K$ so that also each spine is in general position with respect to the sweep-out coming from the other bridge surface and the two sweep-outs themselves are in general position with respect to each other. This operation will be referred to as a two-parameter sweep-out. Two-parameter sweep-outs and their associated graphic are defined in detail in \cite{To} so we only provide a brief overview here.

Associated to a two-parameter sweep-out is a square $I \times I$ where each point $(s,t)$
in the square represents a configuration of $P_K$ and $Q_K$ during the sweep-out. Cerf theory \cite{Ce} describes a graphic $\Gamma$ in the square, a graph all of whose relevant vertices are degree four. The edges of $\Gamma$ correspond to configurations where either
\begin{itemize}
\item$P_K$ and $Q_K$ intersect in a saddle singularity or
\item $P_K, Q_K$ and $K$ intersect transversally at a point.  
\end{itemize}
The former edges are called {\em saddle edges} and the latter (which are best thought of as points where a transverse intersection curve of $P_K$ and $Q_K$ intersects the knot $K$) are called {\em $K$-edges}. Vertices in the graphic correspond to configurations where two of these events occur simultaneously.  That is, either $P_K$ and $Q_K$ intersect in two saddle singularities, or they intersect in a single saddle singularity and elsewhere a curve of $P_K \cap Q_K$ intersects $K$, or $P_K \cap Q_K$ intersects $K$ in two different points. (According to Cerf theory, some edges may also contain two-valent ``birth-death" vertices, but these are irrelevant to our argument.)  Each component of $(I \times I) - \Gamma$, called a region, corresponds to a configuration in which $P_K$ and $Q_K$ are transverse.  

Given a region of a 2-parameter sweep-out of $P_K$ and $Q_K$ as described above, let $\mcCP$ (resp. $\mcCQ$) be the set of all curves of $P_K \cap Q_K$ that are essential in $P_K$ (resp. $Q_K$) and let $\mcC = \mcCP \cup \mcCQ$.  Associate to each region of $(I \times I) - \Gamma$ one or more labels in the following manner:

\begin{itemize}
\item If there are curves in $\mcC$ that are essential on $P_K$ and inessential on $Q_K$ pick an innermost such curve $c$ on $Q_K$ and let $D^c \subset Q$ be the (punctured) disk that $c$ bounds. The interior of $D^c$ might intersect $P_K$, but only in curves that are inessential in both surfaces. If a neighborhood of $\bdd D^c$ lies in $A$ (resp $B$) label the region $A$ (resp $B$).  Analogously label regions with $X$ and/or $Y$.

\item If $\mcC=\emptyset$ (i.e if $P_K$ and $Q_K$ are either disjoint
or only intersect in curves that are inessential on both surfaces), label
the region $a$ if some essential curve on $Q_K$ lies entirely in
$B$. (Note the switch: essential curves of $Q_K$ in $B$ result in the label $a$.)  Use the analogous rule to label regions $b$, $x$ and $y$.

\end{itemize}

\begin{lemma}  \label{lem:onelabel}   If a region in the graphic has label $A$ then there is a curve $c \subset P_K \cap Q_K$ that c-compresses in $A$.  
\end{lemma}

\begin{proof} Let $D^A$ denote a (punctured) disk in $M$ transverse to $P_K$ such that 
\begin{itemize}
\item $\bdd D^A = c^A \subset P_K \cap Q_K$ is essential in $P_K$
\item $interior(D^A)$ intersects $P_K$ only in inessential curves
\item a neighborhood of $\bdd D^A$ in $D^A$ lies in $A$
\item among all such disks with boundary $c^A$, $|D^A \cap P|$ is minimal.
\end{itemize}
Such a (punctured) disk exists by definition of label $A$, in fact it lies in $Q$. Among all the components of $interior(D^A) \cap P$ choose a component $c'$ that is innermost in $P$.  Let $D'^A$ and $D^P$ denote the disjoint (punctured) disks that $c'$ bounds in $D^A$ and $P$ respectively.   Then $S' = D'^A \cup_{c'} D^P$ is a sphere in $M$ with at most two punctures.  $S'$ may not lie completely in any handlebody, but it does intersect $P$ only in inessential curves.  It follows that $S'$ is separating: Any closed curve in $M$ is homologous to a closed curve in $P$, and such a closed curve will intersect $S'$ an even number of times.  Since $S'$ is separating, $|K \cap S'|$ is either $0$ or $2$, so $D'^A$ is punctured if and only if $D^P$ is punctured.  In $D^A$ replace $D'^A$ with $D^P$.  The new (punctured) disk $D'$ still has boundary $c^A$ and all its curves of intersection with $P$ are still inessential in $P_K$ (though curves of $D' \cap Q$ could be essential in $Q_K$),  but $|D' \cap P| < |D^A \cap P|$, since the curve $c'$ and perhaps more have been eliminated.  From this contradiction we deduce that $P$ is disjoint from $interior(D^A)$, ie $D^A \subset A$ as required.  \end{proof}

\begin{lemma}  \label{lem:oneregion}   If a region in the graphic has both labels $A$ and $B$, then $P_K$ is c-strongly compressible.
\end{lemma}

\begin{proof}  Following Lemma \ref{lem:onelabel} there are curves $c^A$ and $c^B$ in $P_K \cap Q_K$ that c-compress in $A$ and $B$ respectively.  Since they arise as curves of intersection of two embedded surfaces, they are disjoint, so the (punctured) disks they bound comprise a c-strong compression of $P$.
\end{proof}

\begin{lemma}  \label{lem:tworegion}   If the union of the labels in adjacent regions in the graphic contains both labels $A$ and $B$, then $P_K$ is c-strongly compressible.
\end{lemma}

\begin{proof}  If both labels are in the same region, the result follows from Lemma \ref{lem:oneregion}.  Otherwise, one region has label $A$ but not $B$ and the adjacent region has label $B$ but not $A$.  

Suppose first that the edge is a saddle edge.  Going through the edge of the graphic separating the regions corresponds to banding two curves $c_+$ and $c_-$ together to give a curve $c$. The three curves cobound a pair of pants on $P_K$ and are thus disjoint. Therefore the curve $c_A$ giving rise to the label $A$ is disjoint from the curve $c_B$ which gives rise to the label $B$. Now the argument in the proof of Lemma \ref{lem:onelabel} shows that these curves bound c-disks for $P_K$ on opposite sides, so $P_K$ is c-strongly compressible.

Now suppose that the edge is a $K$-edge.  Then going through the edge of the graphic corresponds to isotoping a curve of intersection $c$ across a point in $K$.  Beforehand, $c$ gives rise to a label $A$ and after to a label $B$.  But the curves of intersection, before and after the isotopy, are disjoint in $P$, so the same argument applies.
\end{proof}

\begin{lemma}  \label{lem:mixedregion}   If the union of the labels in adjacent regions in the graphic contains both labels $A$ and $b$, then $P_K$ is c-strongly compressible.
\end{lemma}

\begin{proof}  No single region can have both labels, since one implies there is a curve of intersection that is essential in $P_K$ and the other implies that all curves of intersection are inessential in both surfaces.  So suppose one region $R^A$ of the graphic has label $A$ and an adjacent region $R^b$ has label $b$. In the first region there is at least one curve in $\mcC$ that is essential on $P_K$ and inessential on $Q_K$ and in the second all curves of $P_K \cap Q_K$ are inessential on both surfaces. Again the case of a saddle edge is representative and the more difficult, so we describe just that case:

Passing through the edge of the graphic separating the regions corresponds to a saddle tangency where curves $c_+$ and $c_-$ join together to form $c$.  Let $R \subset P_K$, $T \subset Q_K$ be the pairs of pants bounded by the three curves $c_{\pm}, c$ in the two surfaces.  

\medskip

{\bf Case 1: There are no other curves of intersection in the regions $R^A$ and $R^b$.} 

In this case, suppose first that the region in which both $c_+$ and $c_-$ appear is $R^b$ and the region in which $c$ appears is $R^A$.  Then the (punctured) disk $D$ in $Q_K$ bounded by $c$ lies entirely in $A$, since there are no other curves of intersection.  After the saddle move that creates $c_{\pm}$,  $D$ remains as a c-disk for $P_K$ that is disjoint from $Q_K$.  The result then follows from Lemma \ref{lem:disjoint3}.

Next suppose that the region in which both $c_+$ and $c_-$ appear is $R^A$ and the region in which $c$ appears is $R^b$.  With no loss of generality suppose $c_+$ is the curve that gives rise to the label $A$, so $c_+$ is essential on $P_K$ but inessential on $Q_K$.  If the (punctured) disk $D$ that $c_+$ bounds in $Q$ does not contain $c_-$ the proof is just as before.  If it does contain $c_-$ then $c_-$ is also inessential in $Q_K$ and bounds a (punctured) disk in $B$.  If $c_-$ is essential in $P_K$ then the result follows from Lemma \ref{lem:oneregion}.  So suppose $c_-$ is inessential in $P_K$.  Since $c_+$ and $c_-$ are nested in $D$, either $c_-$ or the annulus between $c_+$ and $c_-$ in $D \subset Q_K$ contains no puncture.  In other words, either the (punctured) disk in $Q_K \cap B$ bounded by $c_-$ or the (punctured) disk in $Q_K \cap A$  bounded by $c$ contains no puncture.  It follows that either the (punctured) disk in $P_K$ bounded by $c_-$ or the (punctured) disk in $P_K$  bounded by $c$ (both of which curves we now are assuming are inessential in $P_K$) contains no puncture.  But then $c_+$ also bounds a (punctured) disk in $P_K$, contradicting the fact that it is essential.  

{\bf Case 2:  There are curves of intersection other than $c_{\pm}, c$ in the regions $R^A$ and $R^b$.}

The proof is by induction on the number of such curves. Because of the label $b$ all such curves are inessential in both surfaces; let $\alpha$ be an innermost one in $Q$.  If the (punctured) disk $D^Q$ that $\alpha$ bounds in $Q$ contains $T$ (though, by assumption, no other curves of intersection) the proof is essentially as in Case 1.  If $D^Q$ does not contain $T$, apply the argument of Corollary \ref{cor:reminess} to $\alpha$, reversing the roles of $P$ and $Q$.  The disk $D^P$ that $\alpha$ bounds in $P$ can't contain $R$, since $R$ contains essential curves.  Thus the isotopy of $D^P$ across $D^Q$ that removes $\alpha$ (and perhaps other curves of $D^P \cap Q$) has no effect on the three relevant curves $c_{\pm}, c$.  So $\alpha$ can be removed without affecting the hypotheses.  \end{proof}

\section{Vertices in the graphic with four adjacent labels}

In order to consider labels around vertices in the graphic, return to Conway spheres.  

\begin{defin}  \label{def:tanglesphere} 
A collar $(S, K \cap S) \times I$ of a Conway sphere $S$ is {\em well-placed} with respect to a bridge surface $P$ for $K$ if $P$ intersects the collar as in a Conway product.  That is, up to homeomorphism, $P \cap S \times \{ 0 \}$ is a horizontal great circle $c_0$, $P \cap S \times \{ 1 \}$ is a vertical great circle $c_1$, and $P \cap S \times I$ is a twice-punctured spanning annulus that is the trace of a rotation $\rho_t$ from $c_0$ to $c_1$.  (See Figure  \ref{fig:Conwayprodb}.) 
\end{defin} 

\begin{cor} \label{cor:tanglesphere}  If a bridge surface $P$ is well-placed with respect to a collar $S \times I$ of a Conway sphere $S$, and neither component of $P \cap (S \times \bdd I)$ bounds an unpunctured disk in $P_K$, then either $P_K$ is c-strongly compressible or the surface $P_K - (S \times I)$ is c-incompressible in $M_K - (S \times I)$.  
\end{cor}

\begin{figure}[tbh]
    \centering
    \includegraphics[scale=0.5]{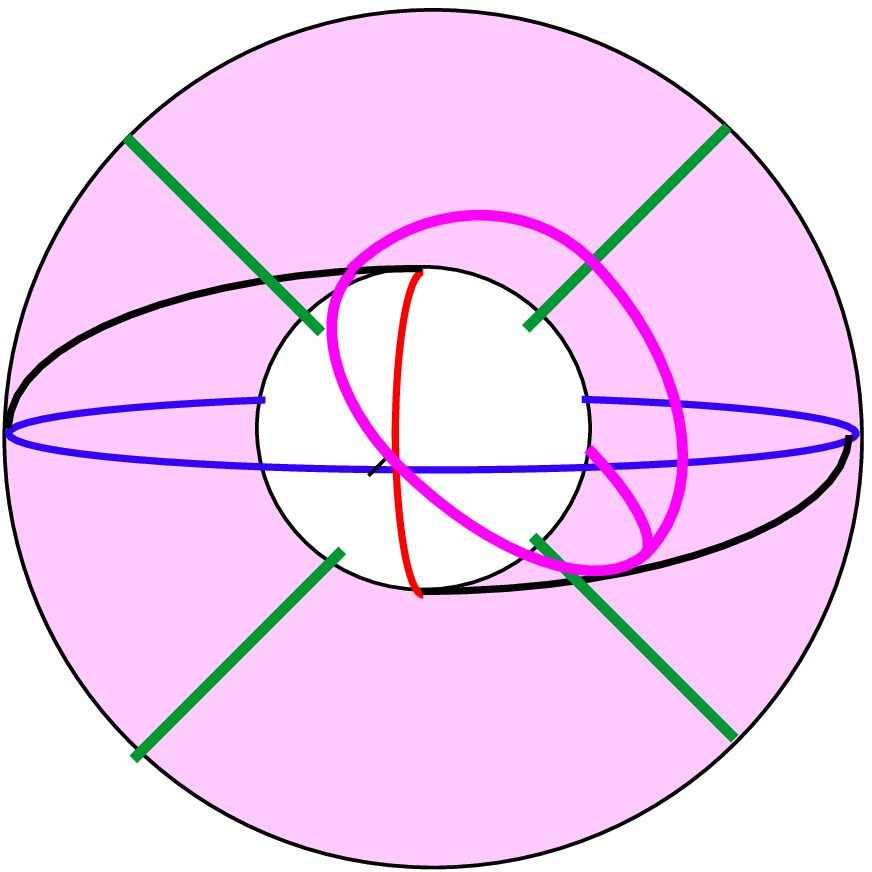}
    \caption{} \label{fig:tanglesphere}
    \end{figure}

\begin{proof} Suppose there is a c-disk for $P_K - (S \times I)$.  The boundary of the c-disk must also be essential in $P_K$, since $P_K$ intersects $S \times I$ in a twice-punctured annulus.  So it suffices to show that there are c-disks in both $A \cap (S \times I)$ and $B \cap (S \times I)$.  Figure \ref{fig:tanglesphere} exhibits a c-disk on one side; one on the other side is symmetrically placed.  (The boundaries of the two disks intersect in two points, so they do not themselves give a c-strong compression.)  
\end{proof}

\begin{lemma}  \label{lem:tanglesphere} Suppose $M = A \cup_P B$ is a bridge decomposition of $K \subset M$ and there is a $4$-punctured sphere $S \subset M$ intersecting $P$ in a single circle $c$ that is essential in $S_K$.  Then either 
\begin{itemize}
\item $c$ bounds a disk in $A_K$ or $B_K$
\item $S_K$ is incompressible (and so is an incompressible Conway sphere)
\item a component of $S_K - c$ is parallel to a component of $P_K - c$ or
\item $P_K$ is c-strongly compressible.
\end{itemize}
\end{lemma}  

\begin{proof}  
$c$ divides $S$ into two twice-punctured disks, $D^A \subset A$ and $D^B \subset B$.  If either are c-ompressible in their respective handlebodies, the c-disk must have boundary parallel to $c$ in the disk, so by parity it must be an unpunctured disk, giving the first conclusion.  

So henceforth we assume that  both $D^A$ and $D^B$ are c-incompressible in their respective handlebodies.  On the other hand,

\medskip

{\bf Claim:  $D^A_K$ (resp $D^B_K$) is $\bdd$-compressible in $A$ (resp $B$)}

Consider how $D^A_K$ intersects a bridge disk $E$ for a bridge that $D^A_K$ intersects.  A circle of intersection that is innermost in $E$ could either be removed by a rechoice of $E$ or the circle is essential in $D^A_K$.  In the latter case, the subdisk of $E$ that the circle bounds would give a compression of $D^A_K$ in $A$, a contradiction.  

So assume there are no circles of intersection and consider an arc of intersection $\alpha$ that is outermost in $E$.  If both ends of the arc lie on $P$, then either the arc can be removed by a rechoice of $E$ or the disk it cuts off from $E$ is a $\bdd$-compressing disk for $D^A_K$ and we are done with the claim.  If both ends of $\alpha$ lie in $K$ then a regular neighborhood of the subdisk of $E$ that $\alpha$ cuts off contains a compressing disk for $D^A_K$ in $A$, a contradiction. If one end of $\alpha$ lies in each of $K$ and $P$, then a regular neighborhood of the subdisk of $E$ that $\alpha$ cuts off contains a $\bdd$-compressing disk for $D^A_K$, establishing the claim. 

\medskip

Either any boundary compression of $D^A_K$ into $A$ produces a c-disk for $P_K$ in $A_K - S$, or $D^A_K$ is parallel in $A_K$ to a twice-punctured subdisk of $P_K$, yielding the third conclusion.  So henceforth assume that any boundary compression of $D^A_K$ into $A_K - S$ (or, symmetrically $D^B_K$ into $B_K - S$) produces a c-disk for $P_K$.  (In particular, following the claim,  $P_K - S$ is c-compressible into  both $A_K - S$  and $B_K - S$.)  If $\bdd$-compressing disks for $D^A_K$ and $D^B_K$ abut $P_K - S$ in disjoint arcs (so, in the terminology of \cite{Sc}, $P_K$ is strongly $\bdd$-compressible to $S_K$), then these boundary compressions create  a pair of c-strong compressing disks for$P$, giving the fourth conclusion.  So to prove the lemma, it suffices to show that if $S_K$ is compressible in $M_K$ then $P_K$ is strongly $\bdd$-compressible to $S_K$.  This we now do:  

Consider the intersection of $P$ with a compressing disk $D$ for $S_K$.  The argument is analogous to the proof of \cite[Theorem 5.4]{Sc} (see also the proof of Lemma \ref{lem:disjoint1} above) so we abbreviate:  Since $D^A_K$ and $D^B_K$ are incompressible in $M_K$, $\bdd D$ must intersect $c$ and $D \cap P$ contains arcs of intersection.  If it contains closed curves of intersection, innermost ones reflect compressing disks for $P_K$ in either $A_K - S$ or $B_K - S$ and if both types arise then $P_K$ is indeed strongly compressible.  If, on the other hand, all innermost disks of $D - P$ cut off compressing disks in $A$, say, regard them as vertices of a graph $\Upsilon \subset D$ in which edges are arcs of intersection with a c-disk $\Ddd^B$ for $B_K - S$.  Either there is an isolated vertex of $\Upsilon$, and so $P_K$ is strongly c-compressible, or there is a non-isolated vertex $v$ incident only to simple edges, in which case the simple edge in $\Upsilon$ that corresponds to an arc of $D \cap \Ddd^B$ that is outermost in $\Ddd^B$ among those incident to $v$, can be used to remove $v$.  So either $P_K$ is strongly compressible or, in this manner, eventually all closed curves of intersection of $P_K$ with $D$ can be removed.

 \begin{figure}[tbh]
    \centering
    \includegraphics[scale=0.5]{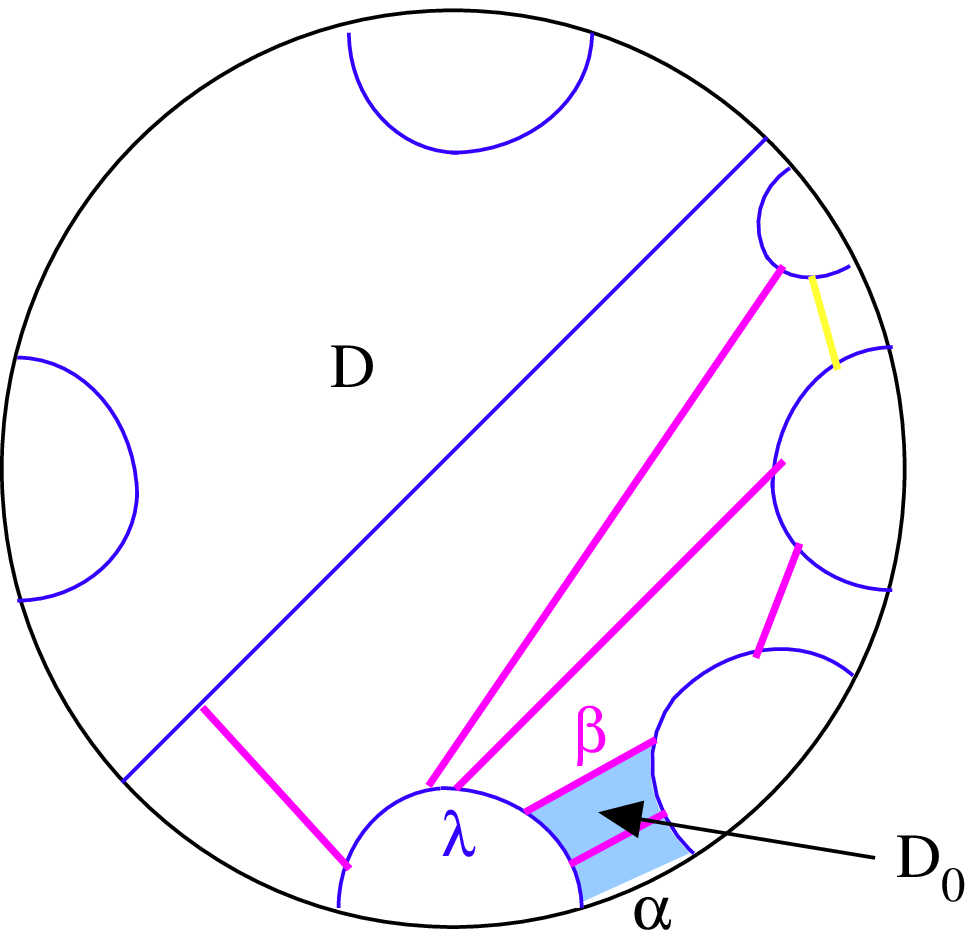}
    \caption{} \label{fig:rectangle}
    \end{figure}

Now consider arcs of intersection.  Suppose with no loss that an arc of $D \cap P_K$ that is outermost in $D$ cuts off from $D$ an outermost disk that lies in $A$.  If there is a similar outermost disk cut off that lies in $B$ then together these describe a strong $\bdd$-compression of $P_K$ to $S_K$ and we are done.  If all outermost arcs of intersection cut off a disk in $A_K$, consider how $\Ddd^B$ intersects $D$.  Among the collection of arcs in $D \cap P$ that are outermost in $D$, there is at least one ($\lambda$ in Figure \ref{fig:rectangle}) that arcs of $D \cap \Ddd^B$ intersects only in arcs that are incident also to adjacent arcs of $P \cap D$.  Let $\bbb$ denote the arc of $D \cap \Ddd^B$ that is outermost in $\Ddd^B$ among all those incident to $\lambda$.  Then there is a rectangle $D_0 \subset (D \cap B)$ with one pair of opposite sides on arcs of $P \cap D$ and other pair of opposite sides consisting of $\bbb$ and an arc $\alpha \subset \bdd D \cap D^B$.  The union of $D_0$ and the disk that $\bbb$ cuts off from $\Ddd^B$ will be a $\bdd$-compressing disk for $P$ to $D^B_K$ that is disjoint from the $\bdd$-compressing disk for $P$ to $D^A_K$ cut off from $D$ by $\lambda$.  This shows that $P_K$ is strongly $\bdd$-compressible to $S_K$, as required.  
\end{proof}

\begin{lemma}  \label{lem:fourlabel}
Suppose there is a $4$-valent vertex in the graphic, the adjacent four regions each have a single label, and, in order around the vertex, these labels are $A$, $X$, $B$, $Y$.  Then either 
\begin{itemize}
\item $P_K$ or $Q_K$ is c-strongly compressible
\item $P_K$ and $Q_K$ can be isotoped to intersect in a non-empty collection of curves that are essential in both surfaces or
\item there is an incompressible Conway sphere $S \subset M$ with a collar $S \times I$ that is well-placed with respect to both $P$ and $Q$.  The pairs of curves $P \cap (S \times \bdd I), (S \times \bdd I)$ are parallel in one component of $S_K \times \bdd I$ and anti-parallel in the other.   The surfaces $P_K - (S \times I)$ and $Q_K - (S \times I)$ are each c-incompressible in $M_K - (S \times I)$.  $P$ and $Q$ can be properly isotoped outside the collar so that all curves of intersection outside the collar are essential in both $P_K$ and $Q_K$.  
\end{itemize}
\end{lemma}

\begin{proof}  We will suppose $P_K$ and $Q_K$ are both c-weakly incompressible and show that only the third or fourth option is possible.  Since all four labels appear, it follows from Lemma \ref{lem:unknot} that $K$ is not the unknot in $S^3$.  The case where one or more edges is a $K$-edge is fairly easy to rule out, so again we focus on the case in which each edge adjacent to the $4$-valent vertex in the graphic corresponds to a saddle move, i.e. attaching a band to a curve or a pair of curves. Let $c^A$, $c^X$, $c^B$ and $c^Y$ denote the collection of all curves of $P_K \cap Q_K$ in each of the corresponding regions. The proof of Lemma \ref{lem:oneregion} shows that no curve in $c^A$ responsible for the label $A$ can be isotoped to be disjoint from a curve in $c^B$ responsible for the label $B$.  It follows that the two disjoint bands by which  $c^A$ is converted to $c^B$ must be incident to the opposite side in $P_K$ of some curve $c_0$ in $c^A$, a curve that is essential in $P_K$ and inessential in $Q_K$.  Let $b^X$ (resp $b^Y$) be the band attached to $c^A$ that changes $c^A$ to  $c^X$ (resp $c^Y$).  

Consider how many curves in $c^A$ are involved in the band moves at the vertex.  There are two bands and each has two ends.  Since $c_0$ is incident to two of those ends, at most two other curves are involved. 

{\bf Case 1:} $c_0$ is the only curve in $c^A$ that is involved in the band moves and also only one curve in $c^B$ is involved in the band moves.  

This means that both ends of $b^X$ and both ends of $b^Y$ are attached to $c_0$.  If the ends of the two bands are not interleaved in $c_0$ (that is, if the ends of $b^X$ lie in a subinterval of $c_0$ that is disjoint from the ends of $b^Y$) then the result of attaching both bands would be three curves in $c^B$, contrary to the assumption in this case.  So the ends of $b^X$ and $b^Y$are interleaved in $c_0$.  Viewed in $Q_K$, where $c_0$ bounds a (punctured) disk, the band moves appear as follows: Two bands are attached to the (punctured) disk of $A \cap Q$ bounded by $c_0$, with the ends of the bands interleaved.  Moreover, the result of attaching both of them is an inessential curve bounding a (once) punctured curve in $B$.  This is only possible if $Q$ is a twice-punctured torus and the bands correspond to orthogonal curves in the torus (call them the meridian and longitude).  The labels $X$ and $Y$ show that the meridian and longitude each bound (punctured) disks on one side or the other, so we further conclude that $Q$ is an unknotted torus in $S^3$.  Since each of the four quadrants has only one label, every curve of intersection that is not incident to either band is inessential in both $P_K$ and $Q_K$.  In particular, in the region labelled $A$, all curves of intersection other than $c_0$ can be removed by a proper isotopy of $P$, following Corollary \ref{cor:reminess}.  

Now do the band move that creates the pair of curves $c^X$, two parallel essential curves in $Q$ with a single puncture in each annulus between them in the torus $Q$.  At that point, consider the surfaces $P^X = P \cap X$ and $P^Y = P \cap Y$.  Exploiting the bridge positioning of $K$ in $X$ and $Y$, let $D^X$, $D^Y$ be meridian disks for the solid tori $X$ and $Y$ respectively, chosen to be disjoint from $K$ and, subject to that, have minimal intersection with $P^X$ and $P^Y$ respectively.  

Since the region is labeled $X$, one component of $P^X$ is a (punctured) disk; let $P^0$ be the other component.  A standard cut and paste argument shows that $D^X$ is disjoint from the (punctured) disk component of $P^X$.  It follows that $\bdd D^X$ is parallel to $\bdd P^0$ in $Q_K$ so either it reveals a compressing disk for $P^0_K$ in $X_K$ or $P^0_K$ is itself an unpunctured disk.  But in the latter case, the original $c_0$ would  bound the band sum of a punctured and an unpunctured disk, ie an unpunctured disk and so it would be inessential in $P_K$, and could not have given the label $A$.  So we conclude that $P^X_K$ is compressible in $X_K$.  

Similarly, $P^Y_K$ is at least c-compressible in $Y_K$. (Recall that we are still considering the region with label $X$, so there is not necessarily any symmetry with $P^X_K$.)   For consider the intersection of $P^Y$ with $D^Y$.   
 If there are any closed curves in $P^Y \cap D^Y$, an innermost one in $D^Y$ would give a compressing disk for $P^Y_K$ in $Y_K$, as we are seeking.  If there are only arcs of intersection, $\bdd$-compression of $P^Y_K$ via an outermost disk cut off in $D^Y$ would reveal either a c-disk for $P^Y_K$ in $Y_K$ or that $P^Y$ is a $\bdd$-parallel punctured annulus.  The latter is impossible: one of the bands creating $c^A$ or $c^B$ from $\bdd P^Y$ must lie in $Y$, else $c^Y$ could not bound a (punctured) disk in $Y$, and that band move would just be the boundary compression, creating only a curve that is inessential in both $P$ and $Q$.

Having established that $P^X_K$ is compressible in $X_K$ and $P^Y_K$ is at least c-compressible in $Y_K$, consider where these c-disks for $P^X_K$ and $P^Y_K$ must lie.  If one lies in $X_K \cap A$ and the other on $Y_K \cap B$, they would give a c-strong compression of $P_K$.  So they both lie on the same side of $P$, say in $X_K \cap A$ and $Y_K \cap A$.  A single band move changes $c^X$ to $c^B$.  If that band is in $X_K$ then the c-disk in $Y_K \cap A$ is undisturbed.  Similarly if that band is in $Y_K$ then the c-disk in $X_K \cap A$ is undisturbed.  In any case, there is a c-disk for $P_K$ in $A$ that is disjoint from $c^B$, contradicting c-weak incompressibility of $P$.  

{\bf Case 2:} Exactly two curves are involved in the band move at $c^A$ and each band $b^X$ and $b^Y$ has one end on each curve.  

Then the ends of the bands $b^X$ and $b^Y$ that are not in $c_0$ both lie in a curve $c_1$.   But this implies that  the only curve in $c^X$ involved in the band move is the curve obtained by banding $c_0$ and $c_1$ together along $b^X$.  Similarly only one curve in $c^Y$ is involved in the band moves.   So this is exactly Case 1, with the roles of $P$ and $Q$ reversed.  

\medskip

{\bf Claim 1: All other cases are equivalent to this one:  exactly three curves of $c^A$ are involved in the band move.}

If instead $c_0$ is the only curve of $c^A$ involved in the band move and, unlike Case 1, the ends of the bands are not interleaved, then the result is that exactly three curves of $c^B$ are involved in the band move.  Switching $A$ and $B$ in the argument establishes the claim in this case.  If exactly one other curve $c_1$ is involved in the band move and, unlike Case 2, the end of one band (say $b^Y$) has both ends on $c_0$, then exactly three curves of $c^Y$ are involved in the band moves.  Thus switching the roles of $P$ and $Q$ and of $A$ and $Y$ establishes the claim also in this case.  

\medskip

Having established the claim, we now proceed under the assumption that exactly three curves of $c^A$ are involved in the band move. Then $b^X$ attaches $c_0$ to another curve $c_{\ell}$.  In order that this band move create the new label $X$, the band $b^X$ must span a (punctured) annulus in $P$ between $c_0$ and $c_{\ell}$.  The same is true for $b^Y$, so all three curves $c_{\ell}, c_0$ and $c_r$ in $c^A$ that are involved in the band moves are co-annular in $P$, $c_0$ is in the middle, and the entire annulus contains at most two punctures.  See Figure \ref{fig:threeinarow}.

\begin{figure}[tbh]
    \centering
    \includegraphics[scale=0.6]{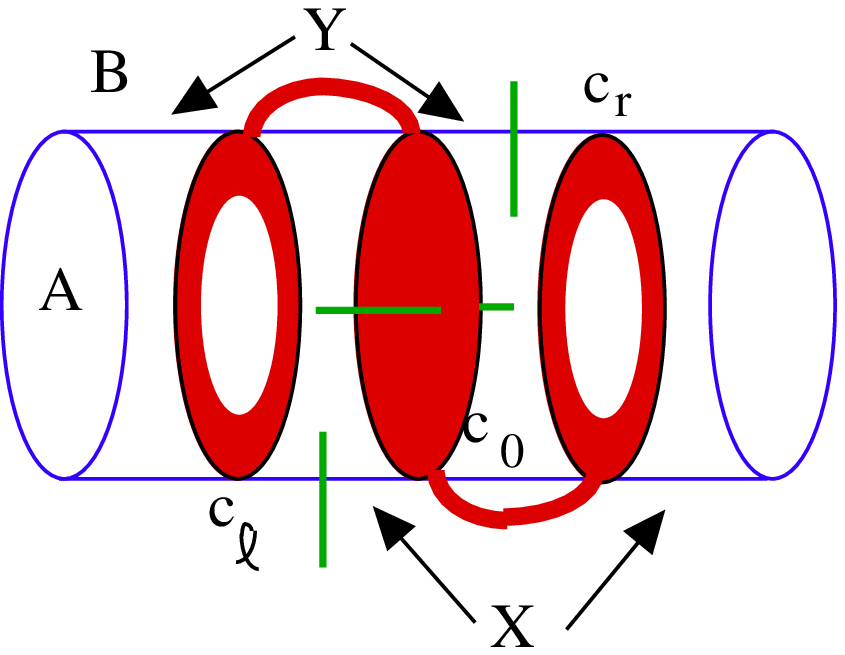}
    \caption{} \label{fig:threeinarow}
    \end{figure}

\medskip

{\bf Claim 2: Both $c_{\ell}$ and $c_r$ are essential in both surfaces. }

If, say $c_r$, were essential in exactly one surface, it would give a label (and so by hypothesis give the label $A$) and that label would persist when $b^X$ is attached, since attaching the band has no effect on $c_r$.  This would contradict the hypothesis that there is only one label in each region.  Suppose, on the other hand, $c_r$ were inessential in both surfaces.  If it bounded an unpunctured disk in either surface, hence in both surfaces, then the band move on $b^Y$ could have no effect; the resulting curve in $P \cap Q$ would be isotopic to $c_0$. Suppose, on the other hand, that $c_r$ bounded a punctured disk in both surfaces.  Then, by parity, the (punctured) annulus between $c_0$ and $c_r$ and the (punctured) disk in $Q$ bounded by $c_0$ must together have exactly one puncture. If the annulus had no puncture, then $c_0$ would, against hypothesis, be inessential in $P$.  If the disk in $Q$ bounded by $c_0$ contained no puncture then when $c_0$ is banded to $c_r$ the result would be, against hypothesis, inessential in $Q$.  The contradictions establish the claim.  

\medskip

{\bf Claim 3: All curves of intersection that are inessential in both $P_K$ and $Q_K$ can be removed by isotopies of $P$ and $Q$ that do not affect $c_0$, $c_{\ell}$ or $c_r$ or either band $b^X, b^Y$}

The proof is by induction on the number of curves inessential in both surfaces. Let $\alpha$ be an innermost one in $Q$.  The (punctured) disk $D$ that $\alpha$ bounds in $Q_K$ can't contain $c_{\ell}$ or $c_r$, since these are essential in $Q$.  Nor can it contain $c_0$ since $c_0$ is banded via the bands $b^X, b^Y$ to the essential curves $c_{\ell}$ and $c_r$ in the complement of $\alpha$.  Apply the argument of Corollary \ref{cor:reminess} to $\alpha$, reversing the roles of $P$ and $Q$.  The disk $D^P$ that $\alpha$ bounds in $P$ can't contain $c_0$, $c_{\ell}$ or $c_r$, nor the bands $b^X, b^Y$ between them, since the curves are essential in $P_K$.  Thus the isotopy of $D^P$ across $D^Q$ that removes $\alpha$ (and perhaps other curves of $D^P \cap Q$) has no effect on the relevant curves and bands.  This isotopy removes $\alpha$ and completes the inductive step.  

\medskip

{\bf Claim 4: Neither $c_{\ell}$ nor $c_r$ can bound  a (punctured) disk in any of $X, Y, A$ or $B$.}

Suppose, say, that $c_{\ell}$ bounded a (punctured) disk in $X$ or $Y$.  Note that $c_{\ell}$ can be isotoped in $Q_K$ to be  disjoint from the curve obtained by banding $c_0$ to either $c_{\ell}$ or $c_r$, so, following Claim 3, it is disjoint from curves bounding punctured disks in $Y$ and $X$.  It would follow that $Q_K$ is strongly compressible,  a contradiction.  

To show that $c_{\ell}$ can't bound a (punctured) disk in $A$ or $B$, note that $c_{\ell}$ can be isotoped in $P_K$ to be disjoint  from $c_0$ and from the curve obtained by banding $c_0$ simultaneously to both $c_{\ell}$ and $c_r$.  Then the same argument applies.  

\medskip

Following the claims, we consider the appearance of the three curves in $Q$:  $c_0$ bounds a (punctured) disk that lies in $Q \cap A$.  When $c_0$ is banded to essential curve $c_r$ the result is still essential in $Q$ (though now inessential in $P$) but when this curve is further banded to $c_{\ell}$ the resulting single curve bounds a (punctured) disk that lies in $Q \cap B$.  Hence $c_r$ and $c_{\ell}$ are also coannular in $Q$, with the annulus between them containing both $c_0$ and at most two punctures, at most one of which is outside the (punctured) disk bounded by $c_0$.  

In fact, all possible punctures in the descriptions above must appear.  If, for example, either the disk bounded by $c_0$ in $Q$ or the annulus bounded by $c_r$ and $c_0$ in $P$ contains no puncture, then $c_r$ also bounds a (punctured) disk in $A_K$, namely the disk obtained by combining the annulus in $P$ and the disk in $Q$ bounded by $c_0$.  But such a disk contradicts Claim 4.  

The next step is to either rearrange the intersection so that $P_K$ and $Q_K$ intersect only in curves essential in both surfaces, conclusion 2, or exhibit the incompressible Conway sphere $S$ of conclusion 3.  Consider the most unusual case first

\medskip

{\bf Claim 5:  If $c_{\ell}$ bounds twice-punctured disks $D^P_K$ and $D^Q_K$ in $P_K$ and $Q_K$ respectively and these twice-punctured disks are parallel in $A \cap Y$, then conclusion 2 holds.}

To demonstrate the claim, use the parallelism to do a $K$-isotopy of $c_{\ell}$ across one of the punctures in $D^P_K$ and $D^Q_K$, adding an extra puncture in the annulus between $c_{\ell}$ and $c_r$ in both $P$ and $Q$.  Then do the band move along $b^X$.  See Figure \ref{fig:Kisotopy}.  The result is to make a curve of intersection that is essential in $P_K$, instead of one that gives rise to the label $X$, since the disks the curve bounds in $P_K$ and $Q_K$ are now both twice-punctured disks.  By hypothesis, all other curves of intersection are also essential in both surfaces, yielding conclusion 2.  

\begin{figure}[tbh]
    \centering
    \includegraphics[scale=0.6]{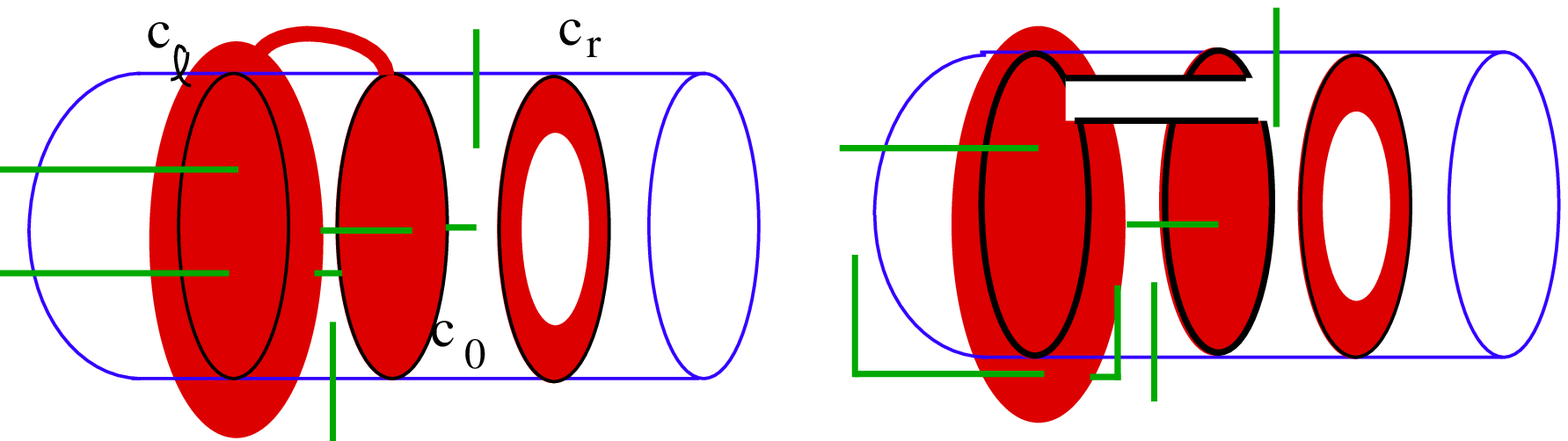}
    \caption{} \label{fig:Kisotopy}
    \end{figure}

A symmetric argument applies if $c_r$ bounds a twice-punctured disk in both $P_K$ and $Q_K$ and these two disks are parallel.  

\medskip

We now construct the relevant Conway sphere.  Start from the description above of how $P$ and $Q$ intersect in the region labeled $A$.  The band $b^X$ must lie in $X$ since a (punctured) disk in $X$ is obtained by removing the band from the annulus in $P$ between $c_0$ and $c_{\ell}$.  Similarly $b^Y$ lies in $Y$.  Push the bands $b^X$ and $b^Y$ close to $Q$ so both lie in a collar $C(Q)$ of $Q$. The end of the collar $Q^X$ in $X$ intersects $P$ as if the band move had been done on $b^X$; the end of the collar $Q^Y$ in $Y$ intersects $P$ as if the band move had been done on $b^Y$.  

 Let $V \subset Q$ be the twice-punctured annulus bounded by $c_{\ell} \cup c_r$; $V$ contains $c_0$ and the once-punctured disk in $Q \cap A$ that it bounds. Choose a curve $c$ in $Q$ just outside of $V$ that is parallel in $Q_K$ to $c_{\ell}$.  Use the collar structure to push $c$ to the copy $Q^X$ of $Q$. Since in $Q^X$ the curve $c_0$ has been banded to $c_{\ell}$, the copy of $c$ in $Q^X$ can be isotoped, intersecting $K$ once but disjoint from $P$, to a curve that bounds a once-punctured disk in $X$ disjoint from $P$, namely the disk that gives rise to the label $X$.  Thus $c$ bounds a twice-punctured disk $D^X$ in $X$ that is disjoint from $P$.  On the other hand, a parallel copy $c'$ of $c_{\ell}$ just {\em inside} $V$, when pushed to the copy $Q^Y$ of $Q$ can be isotoped in $Q^Y$ (crossing $K$ once but disjoint from $P$) to the curve obtained when $b^Y$ connects $c_0$ to $c_r$.  It follows that $c'$ bounds a twice-punctured disk $D^Y$ in $Y$ disjoint from $P$. If the disks $D^X$ and $D^Y$ are attached to the annulus in $Q_K$ between $c$ and $c'$, the result is a $4$-punctured sphere $S_{\ell}$ that intersects both $P$ and $Q$ in a single circle.  The circle is essential in $P_K$ and $Q_K$ since $c_{\ell}$ is.  The same construction can be done at $c_r$ instead of $c_l$.  The resulting sphere $S_r$ is parallel to $S_{\ell}$ in $M_K$ and the construction explicitly shows that the region between them is well-placed with respect to both $P$ and $Q$.  See Figure \ref{fig:fourlabel}.  
  
 \begin{figure}[tbh]
    \centering
    \includegraphics[scale=0.6]{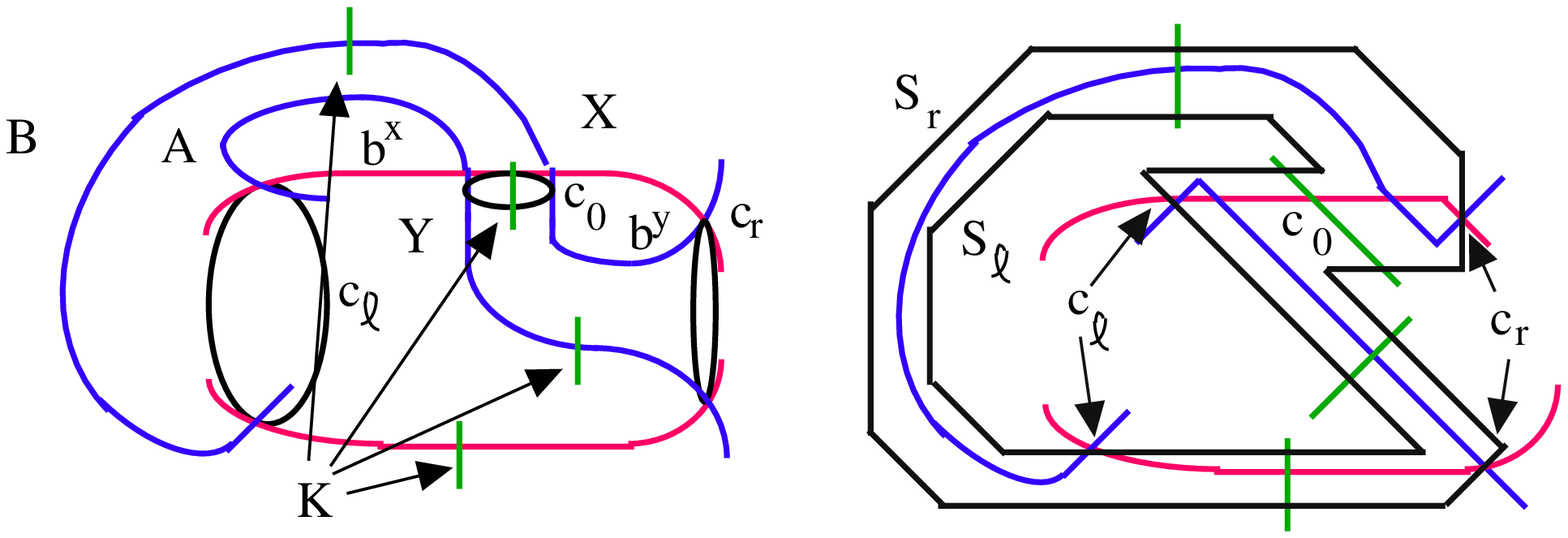}
    \caption{} \label{fig:fourlabel}
    \end{figure}

To determine if $S_{\ell}$  (or $S_r$) is an incompressible Conway sphere for $P$ and $Q$, we apply Lemma \ref{lem:tanglesphere}.  There are three cases:  

\medskip

{\bf Case a: No component of $S_{\ell K} - c_{\ell}$ is parallel to a component of $P_K - c_{\ell}$ or $Q_K - c_{\ell}$.}

Then, following Lemma \ref{lem:tanglesphere}, $S_{\ell}$ is a incompressible Conway sphere for both surfaces.  Corollary \ref{cor:tanglesphere} further shows that the surfaces $P_K - (S \times I)$ and $Q_K - (S \times I)$ are each c-incompressible in $M_K - (S \times I)$.   This is the third possible conclusion sought.

\medskip

{\bf Case b:  $c_{\ell}$ cuts off a component $Q^{\ell}$ of $Q_K \cap A$ that is parallel to $S_{\ell} \cap B.$}

See Figure \ref{fig:fourlabelb}.  In this case, the component $P^{\ell}$ of $P_K \cap Y$ cut off by $c_{\ell}$ lies between the two parallel twice-punctured disks.  It's also incompressible in that region, essentially by Corollary \ref{cor:tanglesphere}.  It follows easily that $P^{\ell}$ is also a parallel twice-punctured disk and so, following Claim 5, conclusion 2 holds.

\begin{figure}[tbh]
    \centering
    \includegraphics[scale=0.6]{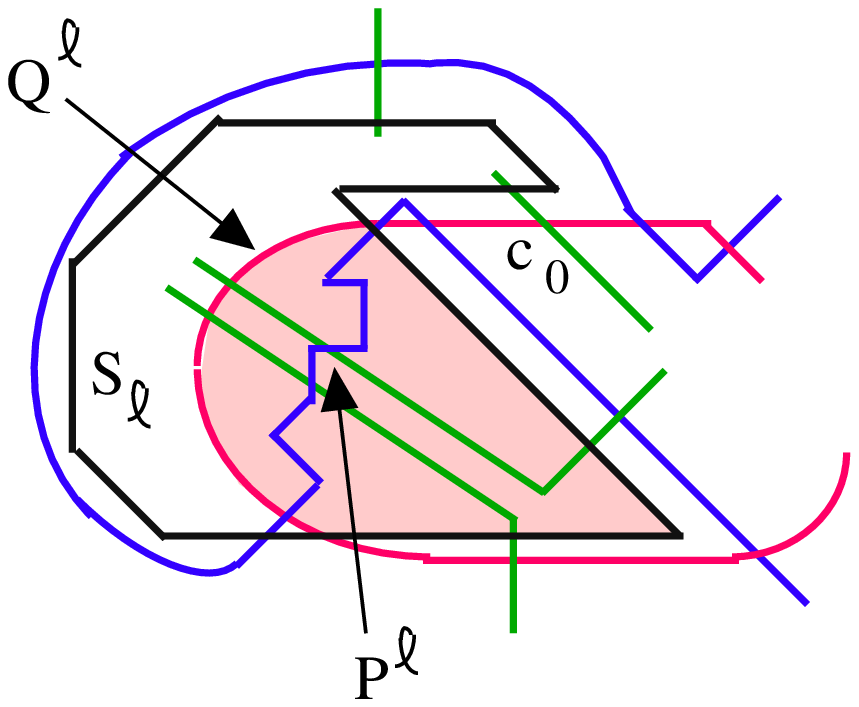}
    \caption{} \label{fig:fourlabelb}
    \end{figure}

The case in which $c_{\ell}$ cuts off a component of $P_K \cap Y$ that is parallel to $S_{\ell} \cap X$ is handled symmetrically.  $c_{\ell}$ can't bound a twice-punctured disk on the other side in either surface (ie the $B$-side of $c_{\ell}$ in $Q_K$ or the $X$ side of $c_{\ell}$ in $P_K$) since $c_r$ is essential in both surfaces.

{\bf Case c:  $c_{\ell}$ cuts off a component $Q^{\ell}$ of $Q_K \cap A$ that is parallel to $S_{\ell} \cap A.$}

See Figure \ref{fig:fourlabelc}. The argument in this case is analogous to that in Claim 5:  The parallelism between the two disks allows us to define a $K$-isotopy that moves the puncture in the annulus between $c_0$ and $c_{\ell}$ in $P$ to the annulus between $c_{\ell}$ and $c_r$ in $Q$.  When both band moves on $b^X$ and $b^Y$ are done (ie in the region labelled $B$ at the vertex) what was previously a punctured disk component of $B \cap Q$ is now a twice-punctured disk, so its boundary is essential in $Q$.  This change, when viewed in $P_K$, simply replaces the old curve of intersection with a curve parallel to $c_r$, which is essential in $P_K$.  Thus conclusion 2 holds again.  

\begin{figure}[tbh]
    \centering
    \includegraphics[scale=0.6]{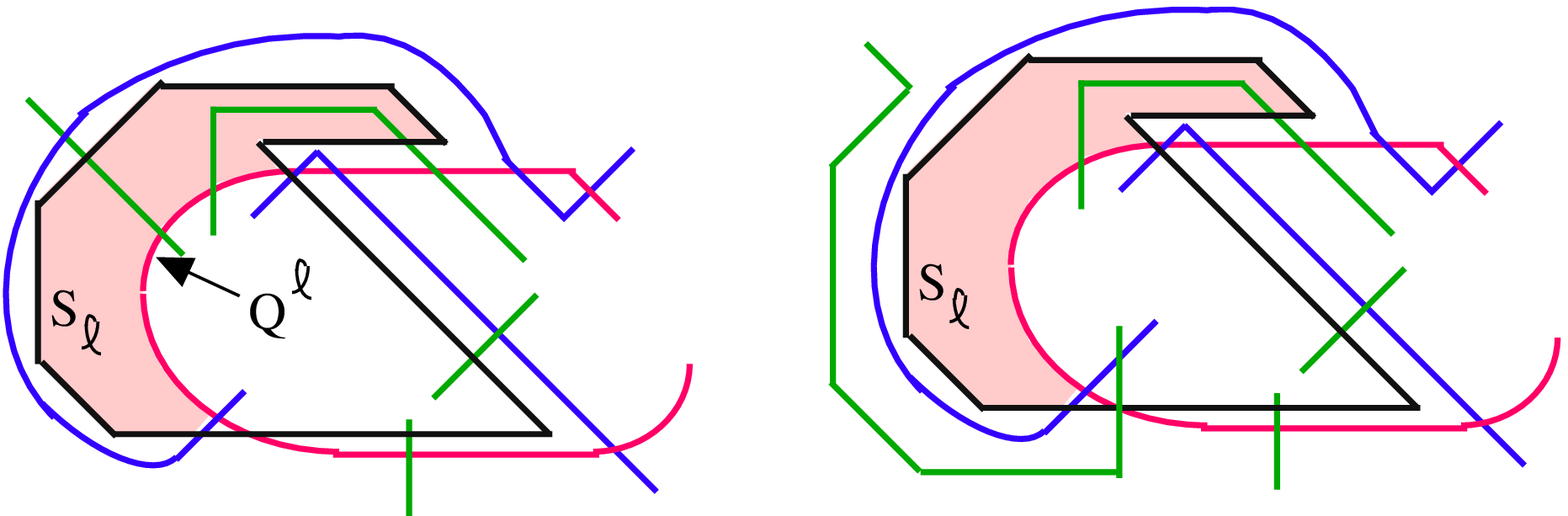}
    \caption{} \label{fig:fourlabelc}
    \end{figure}

The case when $c_{\ell}$ cuts off a component of $P_K \cap Y$ that is parallel to $S_{\ell} \cap Y$ is handled symmetrically.  
\end{proof}

\begin{cor}  \label{cor:fourlabel}
Suppose both $P_K$ and $Q_K$ are c-weakly incompressible and there is a $4$-valent vertex in the graphic, the adjacent four regions each have a single label, and, in order around the vertex, these labels are $A$, $X$, $B$, $Y$.  Then either 
\begin{itemize}
\item $P_K$ and $Q_K$ can be isotoped to intersect in a non-empty collection of curves, each essential in both surfaces or 
\item $P_K$ and $Q_K$ have the same Conway decomposing sphere $S$.  The curves in which the surfaces intersect a collar of $S$ are parallel at one end of the collar and anti-parallel at the other.  Furthermore, the surfaces can be isotoped outside the collar so that all curves of intersection outside the collar are essential in both surfaces.  
\end{itemize}
\end{cor}

Now we can prove the main result of the paper:

\begin{thm} \label{thm:essential}
    Suppose $P$ and $Q$ are bridge surfaces for a link $K$ in a closed orientable $3$-manifold $M$. Assume further that $P_K$ and $Q_K$ are both c-weakly incompressible and neither is a twice-punctured sphere. Then $P_K$ may be properly isotoped so that $P_K$ and $Q_K$ intersect in a non-empty collection of curves so that either
    \begin{itemize} 
    \item all curves of intersection are essential in both surfaces or
    \item $P_K$ and $Q_K$ have the same Conway decomposing sphere $S$.  The curves in which the surfaces intersect the collar of $S$ are parallel at one end of the collar and anti-parallel at the other.  Outside the collar all curves of intersection are essential in both surfaces.
    \end{itemize}
\end{thm}

An important special case is an immediate corollary:

\begin{cor} \label{cor:essential}  Suppose $P$ and $Q$ are bridge surfaces for a link $K$ in a closed orientable $3$-manifold $M$. Assume further that $K \subset M$ is not the unknot in $S^3$, that there is no incompressible Conway sphere for $K$, and that both $P_K$ and $Q_K$ are c-weakly incompressible.  Then $P_K$ may be properly isotoped so that $P_K$ and $Q_K$ intersect in a non-empty collection of curves that are essential in both surfaces.
\end{cor}

{\em Proof of Theorem \ref{thm:essential}.}  Since neither $P_K$ nor $Q_K$ is a twice-punctured sphere, the neighborhood of any bridge disk for any bridge surface $P$ contains a compressing disk for $P_K$.  Since $P_K$ is c-weakly incompressible, it follows from Lemma \ref{lem:unknot} that $K \subset M$ is not the unknot in $S^3$.  

Consider a 2-parameter sweep-out of $P_K$ and $Q_K$ together with the labeling scheme for its graphic described above. The theorem reduces to proving that either there is an unlabeled region or there is a $4$-valent vertex as described in Corollary \ref{cor:fourlabel}.  We will prove the existence of such a region or such a vertex in a sequence of claims.  The labeling is symmetric with respect to the surface we are considering so any statement regarding labels $A$ and $B$ say, has an equivalent statement for labels $X$ and $Y$.  

\medskip

\textbf {Claim 1: The union of the labels of two adjacent regions cannot contain both labels $A$ and $B$.}

This follows immediately from Lemma \ref{lem:tworegion}.  

\medskip

\textbf {Claim 2: The union of the labels of two adjacent regions cannot contain both labels $a$ and $b$}

As in the proof of Lemma \ref{lem:tworegion}, going through the edge of the graphic separating the regions corresponds to banding two curves $c_+$ and $c_-$ together to give a curve $c$. As both regions have small labels all curves of $P_K \cap Q_K$ before and after passing through the saddle are inessential on both surfaces.  Thus the three curves $c_{\pm}, c$ bound a pair of pants $F \subset P$, and each curve bounds a (punctured) disk in $P$.  The (punctured) disks bounded by the three curves can't all be disjoint from $F$, else $P$ would be a twice-punctured sphere, so one of the (punctured) disks $D \subset P$  bounded by the three curves contains all of $F$.   

No essential curve of $P_K$ can lie inside of $D$, so at least parts of both the essential curves $c_a \subset P_K \cap B$ and $c_b \subset P_K \cap A$ that give rise to the labels $a$ and $b$ lie outside $D$.   $c_a$ and $c_b$ are also automatically disjoint from all the other punctured disks in $P$ bounded by components of $P \cap Q$ since $c_a$ and $c_b$ are disjoint from $Q$.  But removing these disks from $P_K$ leaves a connected surface, so parts of $c_a$ and $c_b$ lie in the same component of  $P - Q$.  This is impossible, since $Q$ separates $A$ from $B$ in $M$.  

\medskip

 \textbf {Claim 3:  The union of the labels of two adjacent regions cannot contain both labels $A$ and $b$}

This follows immediately from Lemma \ref{lem:mixedregion}.  

\medskip

\textbf {Claim 4: The theorem is true if there is an unlabelled region}

In the corresponding region all curves of intersection are either essential in both surfaces or inessential in both surfaces and some must be essential in both surfaces to avoid a label $x$ or $y$.  Apply Corollary \ref{cor:reminess} to remove all inessential circles of intersection but leave all essential curves of intersection.  This gives the first conclusion of the theorem.  

\medskip

Relabel the regions of the graphic as follows. Assign a region label $\aA$ if $A$ or $a$ is amongst the labels of the region. Similarly replace $B$ and $b$ by $\bB$,  $x$ and $X$ by $\xX$ and $y$ and $Y$ by $\yY$.

{\em By Claims 1-3, labels $\aA$ and $\bB$ never appear as labels of the same region or labels of
adjacent regions. The same holds for labels $\xX$ and $\yY$.}

\medskip

\textbf {Claim 5: The theorem is true if there is a vertex surrounded by regions that have all four labels $\aA$, $\bB$, $\xX$ and $\yY$.}

Suppose such a vertex exists and label the four regions clockwise $R_n$, $R_e$, $R_s$ and $R_w$. Without loss of generality suppose that $R_n$ is labeled $\aA$. Then according to Claims 1-3, $R_s$ must be the region that carries the label $\bB$ and, furthermore, neither $R_e$ nor $R_w$ can contain either the label $\aA$ or the label $\bB$.  If either $R_e$ or $R_w$ is unlabelled then we are done by Case 4, so with no loss assume that $R_e$ carries the label $\xX$.  Then the symmetric argument shows that $R_w$ must carry the label $\yY$.  

All the labels must come from upper case letters:  Suppose, for example, that $R_n$ was labeled $a$.  This implies that all curves of intersection are inessential in both surfaces, so some essential curve of $P$ lies in either $X$ or $Y$.  Hence $R_n$ would also have the label $y$ or $x$.  Then $R_n$ and either $R_w$ or $R_e$ would contradict Claim 2.

Any curve of intersection other than the ones involved in the saddle moves around the vertex must be either essential in both surfaces or inessential in both, else there would be a specific additional label in all four regions around the vertex, implying one region would have both labels $\xX$ and $\yY$ or both labels $\aA$ and $\bB$.  Either case contradicts Claims 1-3.   Since the labels are in fact upper case $A, B, X, Y$, the claim follows from Corollary \ref{cor:fourlabel}.

   \bigskip

 Consider the labeling of the regions adjacent to $\bdd (I \times I)$. In each of these regions one of the surfaces $P_K$ or $Q_K$ is the boundary of a small neighborhood of one of the spines.  Suppose that
 the north and south boundaries of the square correspond to the spines $\Sigma_{(A,K)}$ and $\Sigma_{(B,K)}$ respectively, and the east and west boundaries correspond to $\Sigma_{(X,K)}$ and $\Sigma_{(Y,K)}$ respectively.  By general position we may arrange that all four spines are disjoint in $M$.  
 
   \medskip

\textbf{Claim 6: No region adjacent to the north boundary of the square has label $\bB$.}

In such a region $A_K$ is a thin neighborhood of the spine $\Sigma_{(A,K)}$ and $Q_K$ sweeps from a neighborhood of $\Sigma_{(X,K)}$ to a neighborhood of $\Sigma_{(Y,K)}$.  When $Q_K$ is near $\Sigma_{(X,K)}$ or $\Sigma_{(Y, K)}$ it lies entirely in the complement of $A_K$, so these regions are both labelled $a$.  It follows from Claim 1 that these regions cannot also have label $\bB$.  To understand the labelling in other regions near the spine, apply general position to the sweep of $Q_K$ across $\Sigma_{(A,K)}$.

During the sweep by $Q_K$, $Q$ goes through no tangencies with $K$, but it can go through vertices of $\Sigma_{(A,K)}$ and have tangencies with edges in $\Sigma_{(A,K)}$.  With no loss, assume that all vertices of $\Sigma_{(A,K)}$ are valence $1$ (where the spine attaches to the extrema of $K$) and valence $3$.  

When $Q$ passes through either a tangency point with an edge of $\Sigma_{(A,K)}$ or through a valence $3$-vertex  of $\Sigma_{(A,K)}$, the effect is to change the number of curves of $P_K \cap Q_K$ that are essential in $P_K$ and bound a disk in $A \cap Q_K$.  (This may require two stages, because perhaps a center tangency is involved, but the center tangency merely creates a curve inessential in both surfaces, so it has no effect on the labeling.)  Thus on one region or the other on either side of the corresponding edge(s) in the graphic, and perhaps on both, there is a label $A$.  It follows from Claims 1 and 2 that on neither region can there be a label $\bB$. 

\begin{figure}[tbh]
    \centering
    \includegraphics[scale=0.6]{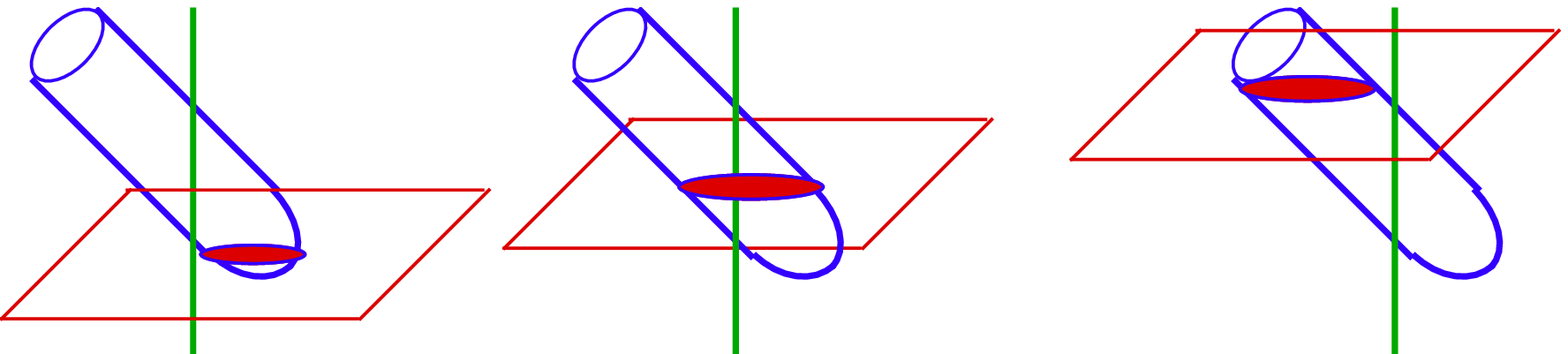}
    \caption{} \label{fig:maxsweep}
    \end{figure}

Suppose $Q$ passes through a valence one vertex of $\Sigma_{(A,K)}$.  The sweep-out of $Q$ is always transverse to $K$, so the sweep-out locally appears as in Figure \ref{fig:maxsweep}.  This is a three-stage process (in which two $K$-edges in the graphic are crossed) that ultimately adds (or removes) the label $A$ but otherwise involves only curves that are inessential in both $P_K$ and $Q_K$.  So if there were also a label $\bB$ in any region involved, there would be one on the region adjacent to label $B$, again contradicting case 1 or 2.

 \bigskip
The proof now follows from the following quadrilateral variant of Sperner's Lemma (\cite[Appendix]{ST}):

 \begin{lemma} \label{lemma:sperner}  Suppose a graph $\Gamma$ is
 properly embedded in a square $I \times I$ and all vertices of $\Gamma$ are
 valence 4 or 1 and all valence one vertices are contained in $\bdd(I
 \times I )$. Suppose further that each region of $I\times I - \Gamma$
 is labeled with letters from the set $\{\aA, \bB, \xX, \yY\}$,
 allowing unlabeled regions and regions with multiple labels, so that

 \begin{itemize}
 \item the union of the labels of two adjacent regions never contains
 both $\aA$ and $\bB$ or both $\xX$ and $\yY$ labels.
 \item  no region adjacent to the north boundary (resp. south, east
 and west boundaries) of $I \times I$ is labeled
 $\bB$ (resp. $\aA$, $\yY$, and $\xX$).
 \end{itemize}

 Then either some region of $(I\times I)-\Gamma$ is unlabeled or all four labels occur in the four regions surrounding some vertex of $\Gamma$. 

 \end{lemma}

\qed

\end{document}